\documentclass[aos,seceqn,citesort,dvips]{arximspdf}
\usepackage{mathbh}
\usepackage{dcolumn}
\usepackage{graphicx}

\doi{10.1214/10-AOS823}
\volume{39}
\issue{1}
\pubyear{2011}
\firstpage{48}
\lastpage{81}

\makeatletter

\newcolumntype{d}[1]{D{.}{.}{#1}}
\newtheorem{lemma}{Lemma}[section]
\newtheorem{theorem}{Theorem}[section]

\makeatother

\begin{document}
\begin{frontmatter}

\title{Regression on manifolds: Estimation of the exterior
derivative\thanksref{T1}}
\runtitle{Regression on manifolds}

\thankstext{T1}{Supported in part by NSF Award CCR-0225610 (ITR),
which supports the CHESS at UC Berkeley, and NSF Grant DMS-06-05236.}

\begin{aug}
\author[A]{\fnms{Anil} \snm{Aswani}\corref{}\ead[label=e1]{aaswani@eecs.berkeley.edu}},
\author[B]{\fnms{Peter} \snm{Bickel}\ead[label=e2]{bickel@stat.berkeley.edu}} and
\author[A]{\fnms{Claire} \snm{Tomlin}\ead[label=e3]{tomlin@eecs.berkeley.edu}}
\runauthor{A. Aswani, P. Bickel and C. Tomlin}
\affiliation{University of California, Berkeley}
\address[A]{Department of Electrical Engineering\\
\quad and Computer Sciences \\
University of California, Berkeley\\
253 Cory Hall\\
Berkeley, California 94720-1770\\
USA\\
\printead{e1}\\
\phantom{E-mail: }\printead*{e3}}
\address[B]{Department of Statistics\\
University of California, Berkeley\\
367 Evans Hall\\
Berkeley, California 94720-3860\\
USA\\
\printead{e2}}
\end{aug}

\received{\smonth{10} \syear{2009}}
\revised{\smonth{4} \syear{2010}}

%
\begin{abstract}
Collinearity and near-collinearity of predictors cause difficulties
when doing regression. In these cases, variable selection becomes
untenable because of mathematical issues concerning the existence and
numerical stability of the regression coefficients, and interpretation
of the coefficients is ambiguous because gradients are not defined.
Using a differential geometric interpretation, in which the regression
coefficients are interpreted as estimates of the exterior derivative of
a function, we develop a new method to do regression in the presence of
collinearities. Our regularization scheme can improve estimation error,
and it can be easily modified to include lasso-type regularization.
These estimators also have simple extensions to the ``large $p$,
small~$n$'' context.
\end{abstract}

%

\end{frontmatter}

\section{Introduction}

Variable selection is an important topic because of its wide set of
applications. Amongst the recent literature, lasso-type regularization
\cite{tibshirani1996,knightfu2000,fanli2001,zou2006} and the Dantzig
selector \cite{candestao2007} have become popular techniques for
variable selection. It is known that these particular tools have
trouble handling collinearity. This has prompted work on extensions
\cite{zouhastie2005}, though further developments are still possible.

Collinearity is a geometric concept: it is equivalent to having
predictors which lie on manifolds of dimension lower than the ambient
space, and it suggests the use of manifold learning to regularize
ill-posed regression problems. The geometrical intuition has not been
fully understood and exploited, though several techniques \cite
{knightfu2000,zouhastie2005,belkinetal2006,bickelli2007} have provided
some insight. Though it is not strictly necessary to learn the manifold
for prediction \cite{bickelli2007}, doing so can improve estimation in
a min--max sense \cite{niyogi2008}.

This paper considers variable selection and coefficient estimation when
the predictors lie on a lower-dimensional manifold, and we focus on the
case where this manifold is nonlinear; the case of a global, linear
manifold is a simple extension, and we include a brief discussion and
numerical results for this case. Prediction of function value on a
nonlinear manifold was first studied in \cite{bickelli2007}, but the
authors did not study estimation of derivatives of the function. We do
not consider the case of global estimation and variable selection on
nonlinear manifolds because \cite{goldbergetal2008} showed that
learning the manifold globally is either poorly defined or
computationally expensive.

\subsection{Overview}

We interpret collinearities in the language of manifolds, and this
provides the two messages of this paper. This interpretation allows us
to develop a new method to do regression in the presence of
collinearities or near-collinearities. This insight also allows us to
provide a novel interpretation of regression coefficients when there is
significant collinearity of the predictors.



On a statistical level, our idea is to learn the manifold formed by the
predictors and then use this to regularize the regression problem. This
form of regularization is informed by the ideas of manifold geometry
and the exterior derivative \cite{misneretal1973,lee2003}. Our idea is
to learn the manifold either locally (in the case of a local, nonlinear
manifold) or globally (in the case of a global, linear manifold). The
regression estimator is posed as a least-squares problem with an
additional term which penalizes for the regression vector lying in
directions perpendicular to the manifold.

Our manifold interpretation provides a new interpretation of the
regression coefficients. The gradient describes how the function
changes as each predictor is changed independently of other predictors.
This is impossible to do when there is collinearity of the predictors,
and the gradient does not exist \cite{spivak1965}. The exterior
derivative of a function \cite{misneretal1973,lee2003} tells us how the
function value changes as a predictor and its collinear terms are
simultaneously changed, and it has applications in control engineering
\cite{sastry1999}, physics \cite{misneretal1973} and mathematics \cite
{spivak1965}. In particular, most of our current work is in
high-dimensional system identification for biological and control
engineering systems \cite{aswani2009,aswani2009b}. We interpret the
regression coefficients in the presence of collinearities as the
exterior derivative of the function.

The exterior derivative interpretation is useful because it says that
the regression coefficients only give derivative information in the
directions parallel to the manifold, and the regression coefficients do
not give any derivative information in the directions perpendicular to
the manifold. If we restrict ourselves to computing regression
coefficients for only the directions parallel to the manifold, then the
regression coefficients are unique and they are uniquely given by the
exterior derivative.

This is not entirely a new interpretation. Similar geometric
interpretations are found in the literature \cite
{massy1965,frank1993,fu2000,knightfu2000,yang2004,fu2008}, but our
interpretation is novel because of two main reasons. The first is that
it is the first time the geometry is interpreted in the manifold
context, and this is important for many application domains. The other
reason is that this interpretation allows us to show that existing
regularization techniques are really estimates of the exterior
derivative, and this has important implications for the interpretation
of estimates calculated by existing techniques. We do not explicitly
show this relationship; rather, we establish a link from our estimator
to both principal components regression (PCR) \cite
{massy1965,frank1993} and ridge regression (RR) \cite
{hoerl1970,knightfu2000}. Links between PCR, RR and other
regularization techniques can be shown \cite
{hocking1976,golub1996,helland1998,fu2000}.

\subsection{Previous work}

Past techniques have recognized the importance of geometric structure
in doing regression. Ordinary least squares (OLS) performs poorly in
the presence of collinearities, and this prompted the development of
regularization schemes. RR \cite{hoerl1970,knightfu2000} provides
proportional shrinkage of the OLS estimator, and elastic net (EN) \cite
{zouhastie2005} combines RR with lasso-type regularization. The
Moore--Penrose pseudoinverse (MP) \cite{knightfu2000} explicitly
considers the manifold. MP works well in the case of a singular design
matrix, but it is known to be highly discontinuous in the presence of
near-collinearity caused by errors-in-variables. PCR \cite
{massy1965,frank1993} and partial least squares (PLS) \cite
{wold1975,frank1993,anderson2009} are popular approaches which
explicitly consider geometric structure.

The existing techniques are for the case of a global, linear manifold,
but these techniques can easily be extended to the case of local,
nonlinear manifolds. The problem can be posed as a weighted, linear
regression problem in which the weights are chosen to localize the
problem \cite{ruppertwand1994}. Variable selection in this context was
studied by RODEO \cite{laffertywasserman2005}, but this tool requires a
heuristic form of regularization which does not explicitly consider
collinearity.

Sparse estimates can simultaneously provide variable selection and
improved estimates, but producing sparse estimates is difficult when
the predictors lie on a manifold. Lasso-type regularization, the
Dantzig selector and the RODEO cannot deal with such situations. The
EN produces sparse estimates, but it does not explicitly consider the
manifold structure of the problem. One aim of this paper is to provide
estimators that can provide sparse estimates when the regression
coefficients are sparse in the original space and the predictors lie on
a manifold.

If the coefficients are sparse in a rotated space, then our estimators
admit extensions which consider rotations of the predictors as another
set of tunable parameters which can be chosen with cross-validation. In
variable selection applications, interpretation of selected variables
is difficult when dealing with rotated spaces, and so we only focus on
sparsity in the original space. Numerical results show that our sparse
estimators without additional rotation parameters do not seem to
significantly worsen estimation when there is no sparsity in the
unrotated space.

The estimators we develop learn the manifold and then use this to
regularize the regression problem. As part of the manifold learning, it
is important to estimate the dimension of the manifold. This can either
be done with dimensionality estimators \cite
{costahero2004,levinabickel2005,heinaudibert2005} or with
resampling-based approaches. Though it is known that cross-validation
performs poorly when used with PCR \cite
{kritchmannadler2008,nadler2008}, we provide numerical examples in
Section \ref{section:numexamp} which show that bootstrapping, to choose
dimension, works well with our estimators. Also, it is worth noting
that our estimators only work for manifolds with integer dimensions,
and our approach cannot deal with fractional dimensions.

Learning the manifold differs in the case of the local, nonlinear
manifold as opposed to the case of the global, linear manifold. In the
local case, we use kernels to localize the estimators which (a) learn
the manifold and (b) do the nonparametric regression. For simplicity,
we use the same bandwidth for both, but we could also use separate
bandwidths. In contrast, the linear case has faster asymptotic
convergence because it does not need localization. We consider a linear
case with errors-in-variables where the noise variance is identifiable
\cite{kritchmannadler2008,johnstonelu2009}, and this distinguishes our
setup from that of other linear regression setups \cite
{tibshirani1996,knightfu2000,fanli2001,zouhastie2005,zou2006}.

\section{Problem setup}
\label{section:setup}
We are interested in prediction and coefficient estimation of a
function which lies on a local, nonlinear manifold. In the basic setup,
we are only concerned with local regression. Consequently, in order to
prove results on the pointwise-convergence of our estimators, we only
need to make assumptions which which hold locally. The number of
predictors is kept fixed. Note that it is possible that the dimension
of the manifold varies at different points in the predictor space; we
do not prohibit such behavior. We cannot do estimation at the points
where the manifold is discontinuous, but we can do estimation at the
remaining points.

Suppose that we would like to estimate the derivative information of
the function about the point $X_0 \in\mathbb{R}^p$, where there are
$p$ predictors. The point $X_0$ is the choice of the user, and varying
this point allows us to compute the derivative information at different
points. Because we do local estimation, it is useful to select small
portions of the predictor-space; we define a ball of radius $R$
centered at $X$ in $p$-dimensions using the notation:
$\mathcal{B}^p_{x, R} = \{v \in\mathbb{R}^p \dvtx\|v-x\|_2^2 <
R\}$.\setcounter{footnote}{1}\footnote{In our notation, we denote
subscripts in lower case. For
instance, the ball surrounding the point $X_0$ is denoted in subscripts
with the lower case $x_0$.}

\begin{figure}

\includegraphics{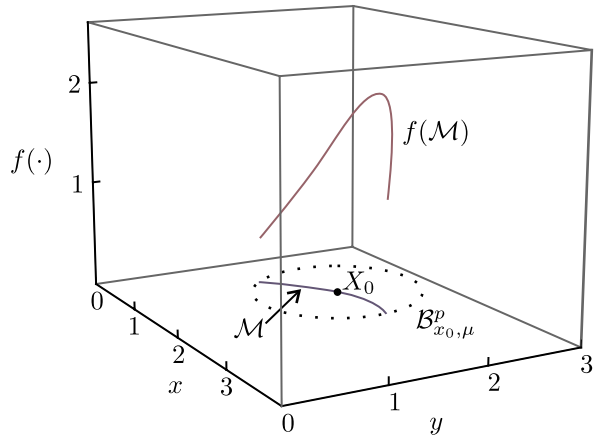}

\caption{In a small ball $\mathcal{B}^p_{x_0,\mu}$ about the
point $X_0$, the predictors form the manifold $\mathcal{M}$. The
response variable is a function of the predictors, which lie on this
manifold. Here the manifold is of dimension $d = 1$, and the number of
predictors is $p = 2$.}
\label{figure:plot}
\end{figure}

We assume that the predictors form a $d$-dimensional manifold $\mathcal
{M}$ in a small region surrounding $X_0$, and we have a function which
lies on this manifold $f(\cdot) \dvtx\mathcal{M} \rightarrow\mathbb{R}$.
Note that $d \leq p$, and that $d$ is in general a function of $X_0$;
however, implicit in our assumptions is that the manifold $\mathcal{M}$
is continuous within the ball. We can more formally define the manifold
at point $X_0$ as the image of a local chart:
%
%
\begin{equation}
\label{eqn:mandef}
\mathcal{M} = \{\phi(u) \in\mathcal{B}^p_{x_0,\mu} \subset\mathbb
{R}^p \dvtx u \in\mathcal{B}^d_{0, r} \subset\mathbb{R}^d\},
\end{equation}
for small $\mu,r > 0$. An example of this setup for $p=2$ and $d=1$ can
be seen in Figure \ref{figure:plot}.

%
%

We make $n$ measurements of the predictors $X_i \in\mathbb{R}^p$, for
$i = \{1, \ldots, n\}$, where the $X_i$ are independent and identically
distributed. We also make $n$ noisy measurements of the function $Y_i =
f(X_i) + \varepsilon_i$, where the $\varepsilon_i$ are independent and
identically distributed with $\mathbb{E}(\varepsilon_i) = 0$ and
$\operatorname{Var}(\varepsilon_i) = \sigma^2$. Let $\kappa, M > 0$ be
finite constants, and assume the following:
\begin{enumerate}
\item The kernel function $K(\cdot)$, which is used to localize our
estimator by selecting points within a small region of predictor-space,
is three-times differentiable and radially symmetric. These imply that
$K(\cdot)$ and $K''(\cdot)$ are even functions, while $K'(\cdot)$ is an
odd function.

\item The bandwidth $h$ is the radius of predictor points about $X_0$
which are used by our estimator, and it has the following asymptotic
rate: $h = \kappa n^{-1/(d+4)}$.

\item The kernel $K(\cdot)$ either has exponential tails or a finite
support \cite{bickelli2007}. Mathematically speaking,
\[
\mathbb{E}\bigl[K^{\gamma}\bigl((X-x)/h\bigr)w(X)\mathbh{1}\bigl(X \in
(\mathcal
{B}^p_{x,h^{1-\varepsilon}})^c\bigr)\bigr] = o(h^{d+4}),
\]
for $\gamma\in\{1, 2\}$, $0 < \varepsilon< 1$ and $|w(x)| \leq M(1
+ |x|^2)$.

\item The local chart $\phi(\cdot)$ which is used to define the
manifold in (\ref{eqn:mandef}) is invertible and three-times
differentiable within its domain. The manifold $\mathcal{M}$ is a
differentiable manifold, and the function $f(\cdot)$ is three-times
differentiable: $\|\partial_i\partial_j\partial_k (f \circ\phi)\|
_{\infty} \leq M$.

\item The probability density cannot be defined in the ambient space
because the Lebesgue measure of a manifold is generally zero. We have
to define the probability density with respect to a $d$-dimensional
measure by inducing the density with the map $\phi(\cdot)$ \cite
{bickelli2007}. We define:
\[
\mathbb{P}(X \in\mathcal{S}) = \mathbb{Q}\bigl(Z \in\phi^{-1}(\mathcal
{S})\bigr),
\]
where $\mathcal{S} \subseteq\mathbb{R}^p$. The density $\mathbb
{Q}(\cdot)$ is denoted by $F(z)$, and we assume that it is three-times
differentiable and strictly positive at $(z = 0) \in\phi^{-1}(X_0)$.

\item The Tikhonov-type regularization parameter $\lambda_n$ is
nondecreasing and satisfies the following rates: $\lambda_n/nh^{d+2}
\rightarrow\infty$ and $h\lambda_n/nh^{d+2} \rightarrow0$. The
lasso-type regularization parameter $\mu_n$ is nonincreasing and
satisfies the following rates: $\mu_n(nh^{d+2})^{-1/2} \rightarrow0$
and $\mu_n(nh^{d+2})^{(\gamma-1)/2} \rightarrow\infty$.
\end{enumerate}

The choice of the local chart $\phi(\cdot)$ is not unique; we could
have chosen a different local chart $\psi(\cdot)$. Fortunately, it can
be shown that our results are invariant under the change of coordinates
$\psi^{-1} \circ\phi$ as long as the measure $\mathbb{Q}(\cdot)$ is
defined to follow suitable compatibility conditions under arbitrary,
smooth coordinate changes. This is important because it tells us that
our results are based on the underlying geometry of the problem.

\section{Change in rank of local covariance estimates}
\label{section:colli}

To localize the regression problem, we use kernels, bandwidth matrices
and weight matrices. We define the scaled kernel $K_h(U) =
h^{-p}K(U/h)$, where $h$ is a bandwidth. The weight matrix centered at
$X_0$ with bandwidth $h$ is given by
\[
W_{x_0} = \operatorname{diag}\bigl(K_h(X_1 - X_0), \ldots, K_h(X_n -
X_0)\bigr),
\]
and the augmented bandwidth matrix is given by $H = H^{1/2}H^{1/2}$, where
\[
H^{1/2} = \sqrt{nh^d}
\left[\matrix{ 1 & 0 \cr0 & h\mathbb{I}_{p \times p} }\right].
\]

If we define the augmented data matrix as
\[
X_{x_0} = \left[\matrix{ 1 & (X_1 - X_0)' \cr
\vdots& \vdots\cr
1 & (X_n - X_0)' }\right],
\]
then the weighted Gram matrix of $X_{x_0}$ is
%
%
\begin{equation}
\label{eqn:weightedgram}
\hat{C}_n \triangleq\left[\matrix{ \hat{C}_n^{11} & \hat{C}_n^{12} \cr
\hat
{C}_n^{21} & \hat{C}_n^{22} }\right] = h^p \cdot
H^{-1/2}X_{x_0}'W_{x_0}X_{x_0}H^{-1/2}.
\end{equation}
A formal statement is given in the \hyperref[app]{Appendix}, but the
weighted Gram
matrix (\ref{eqn:weightedgram}) converges in probability to the
following population parameters:
%
%
\begin{eqnarray}
\label{eqn:cee}
C^{11} & = & F(0)\int_{\mathbb{R}^d}{K(d_u\phi\cdot u)\,du}, \nonumber\\
C^{21} & = & C^{12'} = 0, \\
C^{22} & = & F(0)d_u\phi\cdot\biggl[\int_{\mathbb{R}^d}{K\bigl(d_u\phi
(0) \cdot
u\bigr)uu'\,du}\biggr] \cdot
d_u\phi'.\nonumber
\end{eqnarray}

If we expand the $\hat{C}_n^{22}$ term from the weighted Gram matrix
(\ref{eqn:weightedgram}) into
%
%
\begin{equation}
\label{eqn:sampcov}
\hat{C}_n^{22} = \frac{1}{nh^{d+2-p}}\sum_{i=1}^n{K_h(X_i - X_0)(X_i -
X_0)(X_i - X_0)'},
\end{equation}
it becomes apparent that $\hat{C}_n^{22}$ can be interpreted as a local
covariance matrix. The localizing terms associated with the kernel add
a bias, and this causes problems when doing regression. The
$\hat{C}_n^{11}$ term does not cause problems because it is akin to the
denominator of the Nadaraya--Watson estimator \cite{fangijbels1996}
which does not need any regularizations.

The bias of the local covariance estimate $\hat{C}_n^{22}$ causes
problems in doing regression, because the bias can cause the rank of
$\hat{C}_n^{22}$ to be different than the rank of~$C_n^{22}$. The
change in rank found in the general case of the local, nonlinear
manifold causes problems with MP which is discontinuous when the
covariance matrix changes rank \cite{aitchison1982}. In the special
case of a global, linear manifold, a similar change in rank can happen
because of errors-in-variables. It is worth noting that MP works well
in the case of a singular design matrix.

\section{Manifold regularization}
\label{section:projm}

To compensate for this change in rank, we use a Tikhonov-type
regularization similar to RR and EN. The distinguishing feature of our
estimators is the particular form of the regularizing matrix used. Our
approach is to estimate the tangent plane at $X_0$ of the manifold
$\mathcal{M}$ and then constrain the regression coefficients to lie
close to the principal components of the manifold. The idea for this
type of regularization is informed by the intuition of exterior
derivatives.\footnote{We specifically use the intuition that the
exterior derivative lies in the cotangent space of the manifold, and
this statement can be mathematically written as: $d_xf \in T^*_p\mathcal
{M}$.} An advantage of this regularization is that it makes it easy to
apply lasso-type regularizations, and this combination of the two types
of regularization is similar to EN.

To constrain the regression coefficients to lie close to the manifold,
we pose the problem as a weighted least-squares problem with
Tikhonov-type regularization:
%
%
\begin{equation}
\hat{\beta} = {\arg\min}\|W(Y-X\beta)\|_2^2 + \lambda\|\Pi\beta\|_2^2.
\end{equation}
The matrix $\Pi$ is a projection matrix chosen to penalize $\beta$ for
lying off of the manifold. Contrast this to RR and EN which choose $\Pi
$ to be the identity matrix. Thus, RR and EN do not fully take the
manifold structure of the problem into consideration.

Stated in another way, $\Pi$ is a projection matrix which is chosen to
penalize the components of $\beta$ which are perpendicular to the
manifold. The cost function we are minimizing has the term $\|\Pi\beta\|
_2^2$, and this term is large if $\beta$ has components perpendicular
to the manifold. Components of $\beta$ parallel to the manifold are not
penalized because the projection onto these directions is zero.

Since we do not know the manifold a priori, we must learn the
manifold from the sample local covariance matrix $\hat{C}_n^{22}$. We
do this by looking at the principal components of $\hat{C}_n^{22}$, and
so our estimators are very closely related to PCR. Suppose that we do
an eigenvalue decomposition of $\hat{C}_n^{22}$:
%
%
\begin{equation}
\label{eqn:svdlv}
\hat{C}_n^{22} = \left[\matrix{ \hat{U}^R & \hat{U}^N }\right]
\operatorname
{diag}(\lambda^1, \ldots, \lambda^p)\left[\matrix{ \hat{U}^R & \hat
{U}^N }\right]',
\end{equation}
where $\hat{U}^R \in\mathbb{R}^{p \times d}$, $\hat{U}^N \in\mathbb
{R}^{p \times(p-d)}$ and $\lambda^1 \geq\lambda^2 \geq\cdots\geq
\lambda^p$. Note that\vspace*{1pt} the eigenvalue decomposition always exists
because $\hat{C}_n^{22}$ is symmetric. The estimate of the manifold is
given by the $d$ most relevant principal components, and the remaining
principal components are perpendicular to the estimated manifold.

Because we want the projection matrix $\Pi$ to project $\beta$ onto the
directions perpendicular to the estimated manifold, we define the
following projection matrices
%
%
\begin{eqnarray}
\label{eqn:ptild}
\hat{\Pi} &\triangleq& \hat{U}^N \hat{U}^{N'},\nonumber\\[-8pt]\\[-8pt]
\hat{P} &\triangleq& \operatorname{diag}(0, \hat{\Pi}).\nonumber
\end{eqnarray}
The choice of $d$ is a tunable parameter that is similar to the choice
in PCR. These matrices act as a regularizer because $d$ can always be
chosen to ensure that $\operatorname{rank}(\hat{C}_n^{22} + \lambda
\hat{\Pi}_n) = p$. Furthermore, we have the following theorem
regarding the full regularizing matrix $\hat{P}$:
\begin{theorem}[{[Lemma \ref{theorem:proj}, part (d)]}]
Under the assumptions given in Section~\ref{section:setup}, the
following holds with probability one:
%
%
\begin{equation}
\operatorname{rank}(\hat{C}_n + \lambda_n \hat{P}_n/nh^{d+2}) = p+1.
\end{equation}
\end{theorem}

Our estimators can perform better than PCR because of a subtle
difference. PCR requires that the estimate lies on exactly the first
$d$ most relevant principal components; however, our estimator only
penalizes for deviation from the $d$ most relevant principal
components. This is advantageous because in practice $d$ is not known
exactly and because the principal components used are estimates of the
true principal components. Thus, our regularization is more robust to
errors in the estimates of the principal components. Also, our new
regularization allows us to easily add additional lasso-type
regularization to potentially improve the estimation. PCR cannot be
easily extended to have lasso-type regularization.

We denote the function value at $X_0$ as $f|_{x_0}$. Also, we denote
the exterior derivative of $f(\cdot)$ at $X_0$ as $d_xf|_{x_0}$. Then,
the true regression coefficients are denoted by the vector
%
%
\begin{equation}
\label{eqn:beta}
\beta' = \left[\matrix{ f|_{x_0} & d_xf|_{x_0} }\right].
\end{equation}

The nonparametric exterior derivative estimator (NEDE) is given by
%
%
\begin{equation}
\label{eqn:nede}
\hat{\beta} = \mathop{\arg\min}_{\tilde{\beta}}\{h^p\|W_{x_0}^{1/2}(Y -
X_{x_0}\tilde{\beta})\|_2^2 + \lambda_n\|\hat{P}_n \cdot\tilde{\beta}\|
_2^2\},
\end{equation}
where $\hat{P}_n$ is defined using (\ref{eqn:ptild}) with $\hat{C}_n$.
We can also define a nonparametric adaptive lasso exterior derivative
estimator (NALEDE) as
%
%
\begin{eqnarray}
\label{eqn:nalede}\qquad\quad
\hat{\beta} = \mathop{\arg\min}_{\tilde{\beta}} \Biggl\{h^p\|
W_{x_0}^{1/2}(Y -
X_{x_0}\tilde{\beta})\|_2^2
+ \lambda_n\|\hat{P}_n \cdot\tilde{\beta}\|_2^2 +
\mu_n\sum_{j=1}^p{\frac{1}{\hat{w}_j^\gamma}|\tilde{\beta}_j|}\Biggr\},
\end{eqnarray}
where $\hat{P}_n$ is define in (\ref{eqn:ptild}) using $\hat{C}_n$,
$\hat{w}$ is the solution to (\ref{eqn:nede}), and $\gamma> 0$.

Our estimators have nice statistical properties, as the following
theorem shows.
\begin{theorem}
\label{theorem:nede}
If the assumptions in Section \ref{section:setup} hold, then the NEDE
(\ref{eqn:nede}) and NALEDE (\ref{eqn:nalede}) estimators are
consistent and asymptotically normal:
\[
H^{1/2}(\hat{\beta} - \beta) \stackrel{d}{\rightarrow} C^{\dag}\mathcal
{N}(B', \sigma^2V),
\]
where $B$ and $V$ are, respectively, given in (\ref{eqn:bee}) and (\ref
{eqn:vee}), and $C^{\dag}$ denotes the Moore--Penrose psuedoinverse of
$C$. Furthermore, we asymptotically have that $\hat{\beta}' \in\mathbb
{R} \times T^*_P\mathcal{M}$. The NALEDE (\ref{eqn:nalede}) estimator
has the additional feature that
\[
\mathbb{P}\bigl(\operatorname{sign}(\hat{\beta}) = \operatorname
{sign}(\beta
)\bigr) \rightarrow1.
\]
\end{theorem}

Note that the asymptotic normality is biased because of the bias
typical in nonparametric regression. This bias is seen in both the NEDE
(\ref{eqn:nede}) and NALEDE (\ref{eqn:nalede}) estimators, but
examining $B$ one sees that the bias only exists for the function
estimate $\hat{f}_{x_0}$ and not for the exterior derivative estimate
$d_x\hat{f}|_{x_0}$. This bias occurs because we are choosing $h$ to
converge at the optimal rate. If we were to choose $h$ at a faster
rate, then there would be no asymptotic bias, but the estimates would
converge at a slower rate.

It is worth noting that the rate of convergence in Theorem \ref
{theorem:nede} has an exponential dependence on the dimension of the
manifold $d$, and our particular rates are due to the assumption of the
existence of three derivatives. As is common with local regression, it
is possible to improve the rate of convergence by using local
polynomial regression which assumes the existence of higher-order
derivatives \cite{ruppertwand1994,bickelli2007}. However, the general
form of local polynomial regression on manifolds would require the
choice of a particular chart $\phi(\cdot)$ and domain $\mathcal{U}$.
Local linear regression on manifolds is unique in that one does not
have to pick a chart and domain.

As a last note, recall that the rate of convergence in Theorem \ref
{theorem:nede} depends on the dimension of the manifold $d$ and does
not depend on the dimension $p$ of the ambient space. We might
mistakenly think that this means that the estimator converges in the
``large $p$, small $n$'' scenario, but without additional assumptions
these results are only valid for when $p$ grows more slowly than $n$.
Analogous to other ``large $p$, small $n$'' settings, if we assume
sparsity then we can achieve faster rates of convergence, which is the
subject of the next section.

\section{Large $p$, small $n$}
\label{section:lpsnnpc}
We consider extensions of our estimators to the ``large~$p$, small
$n$'' setting. The key difference in this case is the need to
regularize the covariance matrix. Our NEDE (\ref{eqn:nede}) and NALEDE
(\ref{eqn:nalede}) estimators use the eigenvectors of the sample
covariance matrix, and it is known \cite
{ledoitwolf2003,bickellevina2008} that the sample covariance matrix is
poorly behaved in the ``large $p$, small $n$'' setting. To ensure\vadjust{\goodbreak} the
sample eigenvectors are consistent estimates, we must use some form of
covariance regularization \cite{ledoitwolf2003,zouetal2006,bickellevina2008}.

We use the regularization technique used in \cite{bickellevina2008} for
ease of analysis and because other regularization techniques \cite
{ledoitwolf2003,zouetal2006} do not work when the true covariance
matrix is singular. The scheme in \cite{bickellevina2008} works by
thresholding the covariance matrix, which leads to consistent estimates
as long as the threshold is correctly chosen. We define the
thresholding operator as
\[
T_t(m) = m\mathbh{1}(|m| > t),
\]
and by abuse of notation $T_t(M)$ is $T_t(\cdot)$ applied to each
element of $M$.

The setup and assumptions are nearly identical to those of the fixed
$p$ case described in Section \ref{section:setup}. The primary
differences are that (a) $d,p,n$ increase at different rates toward
infinity, and (b) there is some amount of sparsity in the manifold and
in the function. The population parameters $C_n$, analogous to~(\ref
{eqn:cee}), are functions of $n$ and are defined in nearly the same
manner, except with
$[C_n^{21}]_k = F(0)/2\int_{\mathbb{R}^d}{K(d_u\phi\cdot u)\partial
_{ij}\phi^k u^i u^j\,du}$. Their estimates are now defined
\[
\hat{C}_n = H^{-1}X_{x_0}'W_{x_0}X_{x_0};
\]
compare this to (\ref{eqn:weightedgram}). Just as $C_n$ can be
interpreted as a local covariance matrix, we now define a local
cross-covariance matrix:
\[
R_n = \left[\matrix{R_n^{1} \vspace*{2pt}\cr R_n^{2}}\right] = \left
[\matrix{C_n^{11}\cdot f|_{x_0}
\vspace*{2pt}\cr C_n^{21}\cdot f|_{x_0} + C_n^{22}\cdot d_x
f|_{x_0}}\right],
\]
and the estimates are given by
\[
\hat{R}_n = H^{-1/2}X_{x_0}'W_{x_0}Y.
\]
For the sake of brevity, we summarize the other the differences from
Section \ref{section:setup}. The following things are different:
\begin{enumerate}

\item The manifold $\mathcal{M}_n$, local chart $\phi_n(\cdot)$,
manifold dimension $d_n$, number of predictors $p_n$, and density
function $F_n(\cdot)$ are all functions of $n$. We drop the subscript
$n$ when it is clear from the context. These objects are defined in the
same manner as in Section \ref{section:setup}, and we additionally
assume that the density $F(\cdot)$ is Lipschitz continuous.

\item The asymptotic rates for $d,p,n$ are given by $d = o(\log n)$,
\begin{eqnarray*}
h &=& o\bigl((c_n^{4}n/\log p)^{-1/(4+d)}\bigr),\\
c_n\sqrt{\frac{\log p}{nh^d}} &=& o(1);
\end{eqnarray*}
where $c_n$ is a measure of sparsity that describes the number of
nonzero entries in covariance matrices, exterior derivative, etc.

\item The kernel $K(\cdot)$ has finite support and is Lipschitz
continuous, which implies that
\[
K\biggl(\frac{\phi(hu) - \phi(0)}{h}\biggr) = K(d_u\phi\cdot u) = 0,
\]
for $u \notin\mathcal{B}^{d_n}_{0, \Omega}$. Contrast this to the
second assumption in Section \ref{section:setup}.

\item The local chart $\phi_n(\cdot)$, function $f_n(\cdot)$ and local
(cross-)covariance matrices $C_n,R_n$ satisfy the following sparsity conditions:
%
%
\begin{equation}
\label{eqn:sparselinear}
\sum_{k = 1}^p{\mathbh{1}(Q^k \neq0)} \leq c_n \quad\mbox{and}\quad
|Q^k| \leq M,
\end{equation}
for (derivatives of the local chart; the index $k$ denotes the $k$th
component of the vector-valued $\phi$) $Q^k = \partial_i \phi^k,
\partial_{ij} \phi^k, \partial_{ijm} \phi^k, \partial_{ijmn} \phi^k$;
for (derivatives of the function) $Q^k = [d_xf]_k, \partial_{ik}
(f\circ\phi), \partial_{ijk} (f\circ\phi)$; and for (local covariance
matrices) $Q^k = [C_n]_{ik},[R_n]_{ik}$.

\item The smallest, nonzero singular value of the local covariance
matrix is bounded. That is, there exists $\varepsilon> 0$ such that
%
%
\begin{equation}
\label{eqn:sigbound}
\inf_{n > 0}\Bigl(\inf_{\sigma(\cdot) > 0}\sigma(C_n)\Bigr) >
\varepsilon.
\end{equation}
This condition ensures that the regularized inverse of the local
covariance matrix is well defined in the limit; otherwise we could have
a situation with ever-decreasing nonzero singular values.

\item The Tikhonov-type regularization parameter $\lambda_n$ and the
lasso-type regularization parameter $\mu_n$ have the following
asymptotic rates:
\begin{eqnarray*}
&\displaystyle\lambda_n = O_p\biggl(\sqrt{c_n}\biggl(\frac{nh^d}{\log
p}\biggr)^{1/4}\biggr),&\\
&\displaystyle\mu_n c_n^{3/2}\biggl(\frac{\log p}{nh^d}\biggr)^{1/4} =
o(1), &\\
&\displaystyle\mu_n \biggl(c_n^{3/2}\biggl(\frac{\log p}{nh^d}\biggr
)^{1/4}\biggr)^{1-\gamma} \rightarrow
\infty.&
\end{eqnarray*}

\item The threshold which regularizes the local sample covariance
matrix is given by
%
%
\begin{equation}
t_n = K\sqrt{\frac{\log p}{n}},
\end{equation}
where $\frac{\log p}{n} = o(1)$. This regularization will make the
regression estimator consistent in the ``large $p$, small $n$'' case.
\end{enumerate}

\subsection{Manifold regularization}

The idea is to regularize the local sample covariance matrix by
thresholding. If the true, local covariance matrix is sparse, this
regularization will give consistent estimates. This is formalized by
the following theorem.
\begin{theorem}
\label{theorem:lpsnthresh}
If the assumptions given in Section \ref{section:lpsnnpc} are
satisfied, then
\begin{eqnarray*}
\|T_t(\hat{C}_n) - C_n\| & = & O_p\bigl(c_n\sqrt{\log p/nh^d}\bigr), \\
\|T_t(\hat{R}_n) - R_n\| & = & O_p\bigl(c_n\sqrt{\log p/nh^d}\bigr).
\end{eqnarray*}
\end{theorem}

Once we have consistent estimates of the true, local covariance matrix,
we can simply apply our manifold regularization scheme described in
Section \ref{section:projm}. The nonparametric exterior derivative
estimator for the ``large $p$, small $n$'' case (NEDEP) is given by
%
%
\begin{equation}
\label{eqn:nedep}
\hat{\beta}_n = \mathop{\arg\min}_{\tilde{\beta}} \bigl\|\bigl(T_t(\hat
{C}_n) +
\lambda_n
\hat{P}_n\bigr)\tilde{\beta} - T_t(\hat{R}_n)\bigr\|_2^2,
\end{equation}
where $\hat{P}_n$ is as defined in (\ref{eqn:ptild}) except using
$T_t(\hat{C}_n^{22})$. The nonparametric adaptive lasso exterior
derivative estimator for the ``large $p$, small $n$'' case (NALEDEP) is
given by
%
%
\begin{eqnarray}
\label{eqn:naledep}
\hat{\beta} &=& \mathop{\arg\min}_{\tilde{\beta}} \mathop{\arg\min
}_{\tilde{\beta}}
\bigl\|\bigl(T_t(\hat{C}_n) + \lambda_n \hat{P}_n\bigr)\tilde{\beta} -
T_t(\hat{R}_n)\bigr\|
_2^2 \nonumber\\[-8pt]\\[-8pt]
&&{}+
\mu_n\sum_{j=1}^p{\frac{1}{\hat{w}_j^\gamma}|\tilde{\beta}_j|},\nonumber
\end{eqnarray}
where $\hat{P}_n$ is as defined in (\ref{eqn:ptild}) except using
$T_t(\hat{C}_n^{22})$ and $\hat{w}$ is the solution to (\ref
{eqn:nedep}). These estimators have nice statistical properties.
\begin{theorem}
\label{theorem:nedelpsn}
If the assumptions given in Section \ref{section:lpsnnpc} are
satisfied, then the NALEDE (\ref{eqn:nedep}) and NALEDEP (\ref
{eqn:naledep}) estimators are consistent:
\[
\|\hat{\beta} - \beta\|_2 = O_p\biggl(c_n^{3/2}\biggl(\frac{\log
p}{nh^d}\biggr)^{1/4}\biggr).
\]
The NALEDEP (\ref{eqn:naledep}) estimator is also sign consistent:
\[
\mathbb{P}\bigl(\operatorname{sign}(\hat{\beta}) = \operatorname
{sign}(\beta
)\bigr) \rightarrow1.
\]
\end{theorem}

We do not give a proof of this theorem, because it uses essentially the
same argument as Theorem \ref{theorem:ede}. One minor difference is
that the proof uses our Theorem~\ref{theorem:lpsnthresh} instead of
Theorem 1 from \cite{bickellevina2008}.

\subsection{Linear case}
\label{section:lincase}

Our estimators admit simple extensions in the special case where
predictors lie on a global, linear manifold and the response variable
is a linear function of the predictors. We specifically consider the
errors-in-variables situation with manifold structure in order to
present our formal results, because: in principle, our estimators
provide no improvements in the linear manifold case over existing
methods when there are no errors-in-variables. In practice, our
estimators sometimes provide an improvement in this case. Furthermore,
our estimators provide another solution to the identifiability problem
\cite{fu2008}; the exterior derivative is the unique set of regression
coefficients because the predictors are only sampled in directions
parallel to the manifold, and there is no derivative information about
the response variable in directions perpendicular to the manifold.

Suppose there are $n$ data points and $p$ predictors, and the dimension
of the global, linear manifold is $d$. We assume that $d, n, p$
increase to infinity, and leaving $d$ fixed is a special case of our
results. We consider a linear model $\eta= \Xi\overline{\beta}$,
where $\eta\in\mathbb{R}^{n\times1}$ is a vector of function values,
$\Xi\in\mathbb{R}^{n \times p}$ is a matrix of predictors, and
$\overline{\beta} \in\mathbb{R}^p$ is a vector.

The $\Xi$ are distributed according to the covariance matrix $\Sigma
_{\xi}$, which is also a singular design matrix in this case. The
exterior derivative of this linear function is given by $\beta=
P_{\Sigma_{\xi}}\overline{\beta}$, where $P_{\Sigma_{\xi}}$ is the
projection matrix that projects onto the range space of $\Sigma_{\xi}$.
We make noisy measurements of $\eta$ and $\Xi$:
\begin{eqnarray*}
X &=& \Xi+ \nu,\\
Y &=& \eta+ \varepsilon.
\end{eqnarray*}
The noise $\nu$ and $\varepsilon$ are independent of each other, and
each component of $\nu$ is independent and identically distributed with
mean $0$ and variance $\sigma_{\nu}^2$. Similarly, each component of
$\varepsilon$ is independent and identically distributed with mean $0$
and variance $\sigma^2$. In this setup, the variance $\sigma_{\nu}^2$
is identifiable \cite{kritchmannadler2008,johnstonelu2009}, and an
alternative that works well in practice for low noise situations is to
set this quantity to zero.

Our setup with errors-in-variables is different from the setup of
existing tools \cite{candestao2007,meinshausenyu2006}, but it is
important because in practice, many of the near-collinearities might be
true collinearities that have been perturbed by noise. Several
formulations explicitly introduce noise into the model \cite
{kritchmannadler2008,nadler2008,carrolmacaruppert1999,fantruong1993,johnstonelu2009}.
We choose the setting of \cite{kritchmannadler2008,johnstonelu2009},
because the noise in the predictors is identifiable in this situation.

The exterior derivative estimator for the ``large $p$, small $n$'' case
(EDEP) is given by
%
%
\begin{equation}
\label{eqn:ede}
\hat{\beta} = \mathop{\arg\min}_{\tilde{\beta}} \bigl\|\bigl(T_t(X'X/n)
- \sigma
_{\nu
}^2\mathbb{I} + \lambda_n \hat{P}_n\bigr)\tilde{\beta} - T_t(X'Y/n)
\bigr\|_2^2,
\end{equation}
where $\hat{P}_n$ is as defined in (\ref{eqn:ptild}) except applied to
$\hat{C}_n^{22} = T_t(X'X/n) - \sigma_{\nu}^2\mathbb{I}$. This is
essentially the NEDEP estimator, except the weighting matrix is taken
to be the identity matrix and there are additional terms to deal with
errors-in-variables. We can also define an adaptive lasso version of
our estimator. The adaptive lasso exterior derivative estimator for the
``large $p$, small $n$'' case (ALEDEP) is given by
%
%
\begin{eqnarray}
\label{eqn:alede}
\hat{\beta} &=& \mathop{\arg\min}_{\tilde{\beta}} \bigl\|\bigl
(T_t(X'X/n) - \sigma
_{\nu
}^2\mathbb{I} + \lambda_n \hat{P}_n\bigr)\tilde{\beta} - T_t(X'Y/n)\bigr
\|_2^2
\nonumber\\[-8pt]\\[-8pt]
&&{} + \mu_n\sum_{j=1}^p{\frac{1}{\hat{w}_j^\gamma}|\tilde{\beta
}_j|},\nonumber
\end{eqnarray}
where $\hat{P}_n$ is as defined in (\ref{eqn:ptild}) except applied to
$\hat{C}_n^{22} = T_t(X'X) - \sigma_{\nu}^2\mathbb{I}$ and $\hat{w}$ is
the solution to (\ref{eqn:ede}). We can analogously define the EDE and
ALEDE estimators which are the EDEP and ALEDEP estimators without any
thresholding. In practice, setting the value of the $\sigma_{\nu}^2$
term equal to zero seems to work well with actual data sets.

The technical conditions we make are essentially the same as those for
the case of the local, nonlinear manifold. The primary difference is
that we ask that the conditions in Section \ref{section:lpsnnpc} hold
globally, instead of locally. This also means that we do not use any
kernels to localize the estimators, and the $W$ matrix in the
estimators is simply the identity matrix. If the theoretical rates for
the regularization and threshold parameters are compatibility
redefined, then we can show that these estimators have nice statistical
properties.
\begin{theorem}
\label{theorem:ede}
If the assumptions in Sections \ref{section:lpsnnpc} and
\ref{section:lincase} hold, then the EDEP (\ref{eqn:ede}) and ALEDEP
(\ref
{eqn:alede}) estimators are consistent. They asymptotically converge at
the following rate:
\[
\|\hat{\beta} - \beta\|_2 = O_p\biggl(c_n^{3/2}\biggl(\frac{\log
p}{n}\biggr)^{1/4}\biggr).
\]
The ALEDE (\ref{eqn:alede}) estimator is also sign consistent:
%
%
\begin{equation}
\mathbb{P}\bigl(\operatorname{sign}(\hat{\beta}) = \operatorname
{sign}(\beta
)\bigr) \rightarrow1.
\end{equation}
\end{theorem}

Our theoretical rate of convergence is slower than that of other
techniques \cite{candestao2007,meinshausenyu2006} because we have
applied our technique for local estimation to global estimation, and we
have not fully exploited the setup of the global case. However, we do
get faster rates of convergence in our global case versus our local
case. Furthermore, our model has errors-in-variables, while the model
used in other techniques \cite{candestao2007,meinshausenyu2006} assumes
that the predictors are measured with no noise. Applying the various
techniques to both real and simulated data shows that our estimators
perform comparably to or better than existing techniques. It is not
clear if the rates of convergence for the existing techniques \cite
{candestao2007,meinshausenyu2006} would be slower if there were
errors-in-variables, and this would require additional analysis.

\section{Estimation with data}
\label{section:practical}

Applying our estimators in practice requires careful usage. The NEDE
estimator requires the choice of two tuning parameters, while the
NALEDE and NEDEP estimators require choosing three; the NALEDEP
estimator requires choosing four. The extra tuning parameters---in
comparison to existing techniques like MP or RR---make our method prone
to over-fitting. This makes it crucial to select the parameters using
methods, such as cross-validation or bootstrapping, that protect
against overfitting. It is also important to select from a small number
of different parameter-values to protect against overfitting caused by
issues related to multiple-hypothesis testing \cite{ng1997,juha2003,rao2008}.

Bootstrapping is one good choice for parameter selection with our
estimators \cite
{bickelfreedman1981,shao1993,shao1994,shao1996,bickellevina2008}.
Additionally, we suggest selecting parameters in a sequential manner;
this is to reduce overfitting caused by testing too many models \mbox
{\cite
{ng1997,juha2003,rao2008}}. Another benefit of this approach is that it
simplifies the parameter selection into a set of one-dimensional
parameter searches---greatly reducing the computational complexity of
our estimators. For instance, we first select the
Tikhonov-regularization parameter $\lambda$ for RR. Using the same
$\lambda$ value, we pick the dimension $d$ for the NEDE estimator. The
prior values of $\lambda$ and $d$ are used to pick the
lasso-regu\-larization parameter $\mu$ for the NALEDE estimator.



MATLAB implementations of both related estimators and our estimators
can be found online.\footnote{\url
{http://hybrid.eecs.berkeley.edu/\textasciitilde NEDE/EDE_Code.zip}.} The
lasso-type regressions were computed using the coordinate descent
algorithm \cite{FHHT2007,wulange2008}, and we used the ``improved
kernel PLS by Dayal'' code given in \cite{anderson2009} to do the PLS
regression. The increased computational burden of our estimators, as
compared to existing estimators, is reasonable because of: improved
estimation in some cases, easy parallelization, and computational times
of a few seconds to several minutes on a general desktop for moderate
values of $p$.

\section{Numerical examples}
\label{section:numexamp}

We provide three numerical examples: the first two examples use
simulated data, and the third example uses real data. In the examples
with simulated data, we study the estimation accuracy of various
estimators as the amount of noise and number of data points vary. The
third example uses the Isomap face data\footnote{\url
{http://isomap.stanford.edu/datasets.html}.} used in \cite
{tenenbaumetal2000}. In the example, we do a regression to estimate the
pose of a face based on an image.

For examples involving linear manifolds and functions, we compare our
estimators with popular methods. The exterior derivative is locally
defined, but in the linear case it is identical at every
point---allowing us to do the regression globally. This is in contrast
to the example with a nonlinear manifold and function where we pick a
point to do the regression at. Though MP, PLS, PCR, RR and EN are
typically thought of as global methods, we can use these estimators to
do local, nonparametric estimation by posing the problem as a weighted,
linear regression which can then be solved using either our or existing
estimators. As a note, the MP and OLS estimators are equivalent in the
examples we consider.

Some of the examples involve errors-in-variables, and this suggests
that we should use an estimator that explicitly takes this structure
into account. We compared these methods with Total Least Squares (TLS)
\cite{vanhuffelvandewalle1991} which does exactly this. TLS performed
poorly with both the simulated data and experimental data, and this is
expected because standard TLS is known to perform poorly in the
presence of collinearities \cite{vanhuffelvandewalle1991}. TLS
performed comparably to or worse than OLS/MP, and so the results are
not included.


Based on the numerical examples, it seems that the improvement in
estimation error of our estimators is mainly due to the Tikhonov-type
regularization, with lasso-type regularization providing additional
benefit. Thresholding the covariance matrices did not make significant
improvements, partly because bootstrap has difficulty with picking the
thresholding parameter. Improvements may be possible by refining the
parameter selection method or by changes to the estimator. We also
observed the well-known tendency of lasso to overestimate the number of
nonzero coefficients \cite{mb2009}; using stability selection \cite
{mb2009} to select the lasso parameter would likely lead to better results.

\subsection{Simulated data}
\label{section:lm}

We simulate data for two different models and use this to compare
different estimators. One model is linear, and we do global estimation
in this case. The other model is nonlinear, and hence we do local
estimation in this case. In both models, there are $p$ predictors and
the dimension of the manifold is $d = \operatorname{round}(\frac
{3}{4}p)$. The predictors $\xi$ and response $\eta$ are measured with noise:
\begin{eqnarray*}
x &=& \xi+ \mathcal{N}(0, \sigma_{\nu}^2),\\
y &=& \eta+ \mathcal{N}(0, \sigma^2).
\end{eqnarray*}
And for notational convenience, let $q = \operatorname{round}(\frac
{1}{2}p)$. Define the matrix
\[
F_{ij} = \cases{ 0.3^{|i-j|}, &\quad if $1 \leq i,j \leq d$, \cr
0.3, &\quad if $d+1 \leq i \leq p \wedge j = q+i-d$, \cr
0.3, &\quad if $d+1 \leq i \leq p \wedge j = q+i+1-d$, \cr
0, &\quad o.w.}
\]
The two models are given by:
\begin{enumerate}
\item Linear model: the predictors are distributed $\xi= \mathcal
{N}(0, FF')$, and the function is
%
%
\begin{equation}
\eta= f(\xi) = 1 + \mathop{\sum_{i = 1}}_{i\ \mathrm{is}\ \mathrm
{odd}}^{q}\xi_i.
\end{equation}
If $w = [ 1 \enskip0 \enskip1 \enskip\cdots]$ is a vector with ones
in the odd index-positions and zeros elsewhere, then the exterior
derivative of this linear function at every point on the manifold is
given by the projection of $w$ onto the range space of the matrix~$F$.
\item Nonlinear model: the predictors are distributed $\xi= \sin
(\mathcal{N}(0, FF'))$, and the function is
%
%
\begin{equation}
\eta= f(\xi) = 1 + \mathop{\sum_{i = 1}}_{i \ \mathrm{is}\ \mathrm
{odd} }^{q}\sin(\xi_i).
\end{equation}
We are interested in local regression about the point $x_0 = [0 \enskip
\cdots\enskip0]$. If $w = [1 \enskip0 \enskip1 \enskip\cdots]$ is
a vector with ones in the odd index-positions and zeros elsewhere, then
the exterior derivative of this nonlinear function at the origin is
given by the projection of $w$ onto the range space of the matrix~$F$.
\end{enumerate}

Table \ref{table:lmnoise} shows the average square-loss estimation
error $\|\hat{\beta}-\beta\|_2^2$ for different estimators using data
generated by the linear model and nonlinear model given above, over
%
%
%
\begin{table}[b]
\caption{Averages and standard deviations over 100 replications of
square-loss estimation error for different estimators using data
generated by the linear model and nonlinear model given in Section
\protect\ref{section:lm}, over different noise variances and number of data
points $n$} \label{table:lmnoise}
\begin{tabular*}{\tablewidth}{@{\extracolsep{\fill}}lcccccc@{}}
\hline
& \multicolumn{2}{c}{$\bolds{n = 10}$} & \multicolumn{2}{c}{$\bolds{n = 100}$}
& \multicolumn{2}{c@{}}{$\bolds{n = 1000}$}\\
\hline
& \multicolumn{6}{c@{}}{Linear model: $\sigma_{\nu}^2 = 0.01, \sigma^2 =
1.00$} \\[4pt]
OLS/MP & 3.741 &(2.476)&6.744 &(3.602)& 0.588&(0.368) \\
RR & 2.523 &(1.020)& 0.369 &(0.193)& 0.167&(0.086) \\
EN & 2.562 &(1.054)& 0.117 &(0.197)& 0.017&(0.008) \\
PLS & 2.428 &(0.536)& 0.501 &(0.207)& 0.031&(0.013) \\
PCR & 3.391 &(0.793)& 1.629 &(0.144)& 1.583&(0.047) \\
EDE & 2.528 &(1.030)& 0.367 &(0.185)& 0.166&(0.086) \\
ALEDE & 2.564 &(1.061)& 0.111 &(0.177)& 0.015&(0.006) \\
EDEP & 2.527 &(1.025)& 0.370 &(0.184)& 0.167&(0.085) \\
ALEDEP & 2.563 &(1.057)& 0.111 &(0.177)& 0.015&(0.006)
\\[4pt]
& \multicolumn{6}{c@{}}{Linear model: $\sigma_{\nu}^2 = 0.10, \sigma^2 =
1.00$}\\
[4pt]
OLS/MP & 3.629 &(1.488)& 1.259 &(0.598)& 0.173&(0.050) \\
RR & 2.717 &(0.847)& 0.892 &(0.351)& 0.260&(0.042) \\
EN & 2.783 &(0.895)& 0.661 &(0.523)& 0.064&(0.013) \\
PLS & 2.588 &(0.566)& 0.740 &(0.252)& 0.138&(0.037) \\
PCR & 3.425 &(0.747)& 1.645 &(0.144)& 1.569&(0.047) \\
EDE & 2.716 &(0.846)& 0.840 &(0.305)& 0.255&(0.042) \\
ALEDE & 2.782 &(0.893)& 0.590 &(0.459)& 0.050&(0.013) \\
EDEP & 2.720 &(0.849)& 0.841 &(0.305)& 0.255&(0.042) \\
ALEDEP & 2.784 &(0.895)& 0.590 &(0.459)& 0.050&(0.013) \\
\hline
\end{tabular*}
\end{table}
different noise variances and number of data points $n$. We conducted
100 replications of generating data and doing a regression, and this
helped to provide standard deviations of square-loss\vadjust{\goodbreak} estimation error
to show the variability of the estimators. Table \ref{table:comptime} shows
computation times in seconds for the different estimators; it shows
that our estimators require more computation, but the computation time
is still reasonable. The computation time does not depend on the noise
level, and so we have averaged the computation times over the $100+100$
replications for $\sigma_{\nu}^2=0.01,0.10$.

\setcounter{table}{0}
\begin{table}
\caption{(Continued)}
\begin{tabular*}{\tablewidth}{@{\extracolsep{\fill}}lrrcccc@{}}
\hline
& \multicolumn{2}{c}{$\bolds{n = 20}$} & \multicolumn{2}{c}{$\bolds{n =
100}$} &
\multicolumn{2}{c@{}}{$\bolds{n = 1000}$}\\
\hline
& \multicolumn{6}{c@{}}{Nonlinear model: $\sigma_{\nu}^2 = 0.01, \sigma
^2 = 1.00$}\\[4pt]
OLS/MP & 657.7\phantom{00} &(1842)\phantom{000.}& 3.797 &(2.172)&
0.367&(0.223) \\
RR & 2.188 &(0.715)& 1.160 &(0.265)& 0.300&(0.132) \\
EN & 2.140 &(0.768)& 1.083 &(0.341)& 0.162&(0.081) \\
PLS & 2.102 &(0.757)& 0.995 &(0.256)& 0.172&(0.120) \\
PCR & 2.731 &(0.439)& 1.679 &(0.244)& 0.123&(0.043) \\
NEDE & 2.184 &(0.716)& 1.104 &(0.269)& 0.288&(0.118) \\
NALEDE & 2.136 &(0.769)& 1.009 &(0.357)& 0.144&(0.067) \\
NEDEP & 2.184 &(0.716)& 1.103 &(0.269)& 0.288&(0.118) \\
NALEDEP & 2.136 &(0.769)& 1.008 &(0.357)& 0.144&(0.067)
\\[4pt]
& \multicolumn{6}{c@{}}{Nonlinear model: $\sigma_{\nu}^2 = 0.1, \sigma^2
= 1.00$}\\[4pt]
OLS/MP & 147.5\phantom{00} &(338.3)\phantom{00}& 1.843 &(0.975)&
0.473&(0.110) \\
RR & 2.693 &(1.793)& 1.260 &(0.369)& 0.672&(0.139) \\
EN & 2.698 &(1.927)& 1.168 &(0.455)& 0.472&(0.197) \\
PLS & 2.385 &(0.629)& 1.210 &(0.318)& 0.768&(0.133) \\
PCR & 2.767 &(0.450)& 1.766 &(0.293)& 0.554&(0.348) \\
NEDE & 2.694 &(1.793)& 1.233 &(0.352)& 0.641&(0.119) \\
NALEDE & 2.702 &(1.925)& 1.126 &(0.445)& 0.407&(0.130) \\
NEDEP & 2.693 &(1.794)& 1.231 &(0.352)& 0.641&(0.119) \\
NALEDEP & 2.702 &(1.926)& 1.124 &(0.445)& 0.407&(0.130) \\
\hline
\end{tabular*}
\end{table}

%
%
\begin{table}
\caption{Averages and standard deviations over 200 replications of
computation times in seconds for~the~different~estimators using data
generated by the linear model and
nonlinear~model~given~in~Section~\protect\ref{section:lm}} \label
{table:comptime}
\begin{tabular*}{\tablewidth}{@{\extracolsep{\fill}}lcccccc@{}}
\hline
& \multicolumn{2}{c}{$\bolds{n = 10}$} & \multicolumn{2}{c}{$\bolds{n =
100}$} &
\multicolumn{2}{c@{}}{$\bolds{n = 1000}$}\\
\hline
& \multicolumn{6}{c@{}}{{Linear model}}\\[4pt]
OLS/MP & 0.001 &(0.000)&0.001 &(0.000)& 0.001&(0.000) \\
RR & 0.055 &(0.005)& 0.063 &(0.005)& 0.110&(0.004) \\
EN & 1.739 &(1.325)& 0.292 &(0.032)& 0.387&(0.042) \\
PLS & 1.631 &(0.029)& 1.716 &(0.077)& 1.925&(0.051) \\
PCR & 0.185 &(0.006)& 0.192 &(0.010)& 0.253&(0.007) \\
EDE & 0.239 &(0.007)& 0.255 &(0.016)& 0.367&(0.021) \\
ALEDE & 2.071 &(1.363)& 0.715 &(0.043)& 0.912&(0.059) \\
EDEP & 0.411 &(0.011)& 0.523 &(0.019)& 0.731&(0.052) \\
ALEDEP & 2.248 &(1.366)& 0.990 &(0.042)& 1.283&(0.083) \\
\hline
& \multicolumn{2}{c}{$\bolds{n = 20}$} & \multicolumn{2}{c}{$\bolds{n = 100}$}
& \multicolumn{2}{c@{}}{$\bolds{n = 1000}$}\\
\hline
& \multicolumn{6}{c@{}}{{Nonlinear model}}\\[4pt]
OLS/MP & 0.001 &(0.000)&0.001 &(0.000)& 0.001&(0.000) \\
RR & 0.359 &(0.039)& 0.470 &(0.028)& 0.938&(0.040) \\
EN & 0.916 &(0.433)& 1.005 &(0.058)& 1.917&(0.269) \\
PLS & 1.980 &(0.130)& 2.114 &(0.072)& 2.879&(0.162) \\
PCR & 1.014 &(0.062)& 1.118 &(0.038)& 1.911&(0.166) \\
NEDE & 0.840 &(0.078)& 1.066 &(0.058)& 2.031&(0.092) \\
NALEDE & 2.175 &(0.594)& 1.843 &(0.090)& 3.381&(0.401) \\
NEDEP & 1.403 &(0.120)& 1.726 &(0.087)& 2.277&(0.175) \\
NALEDEP & 2.671 &(0.647)& 2.513 &(0.122)& 4.631&(0.397) \\
\hline
\end{tabular*}
\end{table}

One curious phenomenon observed is that the estimation error goes down
in some cases as the error variance of the predictors $\sigma_{\nu}^2$
increases. To understand why, consider the sample covariance matrix in
the linear case $\hat{S} = X'X/n$ with population parameter $S = FF' +
\sigma_{\nu}^2\mathbb{I}$. Heuristically, the OLS estimate will tend to
$(FF' + \sigma_{\nu}^2\mathbb{I})^{-1}X'Y/n$, and the error in the
predictors actually acts as the Tikhonov-type regularization found in
RR, with lower levels of noise leading to less regularization.

The results indicate that our estimators are not significantly more
variable than existing ones, and our estimators perform competitively
against existing estimators. Though our estimators are closely related
to PCR, RR and EN, our estimators performed comparably to or better
than these estimators. PLS also did quite well, and our estimators did
better than PLS in some cases. Increasing the noise in the predictors
did not seem to significantly affect the qualitative performance of the
estimators, except for OLS as explained above.

In Section \ref{section:lincase}, we discussed how the converge rate of
our linear estimators is of order $n^{-1/4}$ which is in contrast to
the typical convergence rate of $n^{-1/2}$ for lasso-type regression
\cite{meinshausenyu2006}. We believe that this theoretical discrepancy
is because our model has errors-in-variables while the standard model
used in lasso-type regression does not \cite{meinshausenyu2006}. These
theoretical differences do not seem significant in practice. As seen in
Table \ref{table:lmnoise}, our estimators can be competitive with
existing lasso-type regression.

\subsection{Isomap face data}

\begin{figure}

\includegraphics{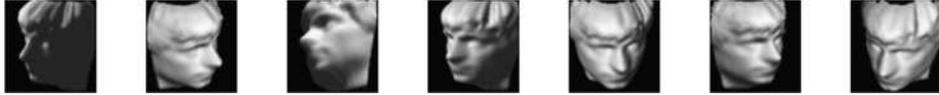}

\caption{Sample images from the Isomap face data
\protect\cite{tenenbaumetal2000}. The images are labeled with and vary
depending upon three variables: illumination-, horizontal- and
vertical-orientation.}
\label{figure:faces}
\end{figure}

The Isomap face data\footnote{\url
{http://isomap.stanford.edu/datasets.html}.} from \cite
{tenenbaumetal2000} consists of images of an artificial face. The
images are labeled with and vary depending upon three variables:
illumination-, horizontal- and vertical-orientation; sample images
taken from this data set can be seen in Figure \ref{figure:faces}.
Three-dimensional images of the face would form a three-dimensional
manifold (each dimension corresponding to a variable), but this data
set consists of two-dimensional projections of the face. Intuitively, a
limited number of additional variables are needed for different views
of the face (e.g., front, profile, etc.). This intuition is confirmed
by dimensionality estimators which estimate that the two-dimensional
images form a low-dimensional manifold \cite{levinabickel2005}.

To compare the different estimators, we did 100 replications of the
following experiment: we randomly split the data ($n = 698$ data
points) into a training set $n_t = 688$ and validation set $n_v = 10$,
and then we used the training set to estimate the horizontal pose angle
of images in the validation set. Since we are doing local linear
estimation, the estimate for each image requires its own regression.
The number of predictors $p$ is large in this case because each data
point $X_i$ is a two-dimensional image. Estimation when $p$ is large is
computationally slow, and so we chose a small validation set size to
ensure that the experiments completed in a reasonable amount of time.
Replicating this experiment 100 times helps to prevent spurious results
due to having a small validation set.

To speed up the computations further, we scaled the images from
$64\times64$ pixels to $7\times7$ pixels in size. This is a
justifiable approach because the images form a low-dimensional
manifold, and so this resizing of the images does not lead to a loss in
predictor information \cite{wright2009}. This leads to significantly
faster computations, because this process reduces the number of
predictors from $p = 4096$ to $p=49$. In practice, our choice of $p=49$
gives predictions that deviate from the true horizontal pose angle of
images (which uniformly ranges between $-$75 to 75 degrees) in the
validation set by a root-mean-squared error of three degrees or less.

%
%

%
%
\begin{table}
\caption{Averages and standard deviations over 100 replications of
validation set prediction error~for~different~estimators using the
Isomap face data \protect\cite{tenenbaumetal2000}. The average computation
time~and~its~standard deviation is given in seconds, and it gives the time
to estimate the~horizontal~pose angle of a single image}
\label{table:pollution}
\begin{tabular*}{\tablewidth}{@{\extracolsep{4in minus
4in}}ld{2.3}d{3.4}d{2.3}c@{}}
\hline
& \multicolumn{2}{c}{\textbf{Prediction error}} &
\multicolumn{2}{c@{}}{\textbf{Computation time}} \\
\hline
OLS/MP & 5.276 &(11.06)&0.002 &(0.001)\\
RR & 5.286 &(6.156)&2.052 &(0.166)\\
EN & 5.168 &(6.112)&17.79 &(5.235)\\
PLS & 10.12 &(14.84)&16.22 &(0.841)\\
PCR & 5.813 &(7.617)&8.877 &(0.609)\\
NEDE & 4.523 &(4.926)&4.740 &(0.347)\\
NALEDE & 4.409 &(4.889)&22.94 &(5.349)\\
NEDEP & 4.527 &(4.925)&7.777 &(0.511)\\
NALEDEP & 4.406 &(4.900)&25.77 &(5.221)\\
\hline
\end{tabular*}
\end{table}

Table \ref{table:pollution} gives the prediction error of the models
generated by different estimators on the validation set. The specific
quantity provided is
%
%
\begin{equation}
{\sum_{X_i \in\mathcal{V}}}\|\hat{\beta}^0(X_i)-Y_i\|_2^2/n_v,
\end{equation}
where $\mathcal{V}$ is the set of predictors in the validation set,
$\hat{\beta}^0(X_i)$ is the first component of the estimated regression
coefficients computed about the point $X_0 = X_i$, and $Y_i$ is the
corresponding horizontal pose angle of the image $X_i$. The regression
is computed using only data taken from the training set. The results
from this real data set shows that our estimators can provide
improvements over existing tools, because our estimators have the
lowest prediction errors. Table \ref{table:pollution} also provides the
computation times for estimating the horizontal pose, and it again
shows that our estimators require more computation but not by an
excessively larger amount.

\section{Conclusion}
By interpreting collinearity as predictors on a lower-dimen\-sional
manifold, we have developed a new regularization, which has connections
to PCR and RR, for linear regression and local linear regression. This
viewpoint also allows us interpret the regression coefficients as
estimates of the exterior derivative. We proved the consistency of our
estimators in both the classical case and the ``large $p$, small $n$''
case and this is useful from a theoretical standpoint.

We provided numerical examples using simulated and real data which show
that our estimators can provide improvements over existing estimators
in estimation and prediction error. Though our estimators provide
modest improvements over existing ones, these improvements are
consistent over the different examples. Specifically, the Tikhonov-type
and lasso-type regularizations provided improvements, and the
thresholding regularization did not provide major improvements. This is
not to say that thresholding is not a good regularization, because as
we showed: from a theoretical standpoint, thresholding does provide
consistency in the ``large $p$, small $n$'' situation. This leaves open
the possibility of future work on how to best select this thresholding
parameter value.

There is additional future work possible on extending our set of
estimators. There is some benefit provided by shrinkage from the
Tikhonov-type regularization which is independent of the manifold
structure. Exploring more fully the relationship between manifold
structure and shrinkage will likely lead to improved estimators.

\begin{appendix}\label{app}
\section*{Appendix}
In this section, we provide the proofs of our theorems. We also give a
few lemmas, which are needed for the proofs, that were not stated in
the main text.
\begin{lemma}
\label{theorem:gram}
If the assumptions in Section \ref{section:setup} hold, then:
\begin{itemize}
\item[(a)] $ \|\mathbb{E}(\hat{C}_n^{22}) - C^{22}\|_2^2 =
O(h^4)$;\vspace*{2pt}
\item[(b)] $ \|\mathbb{E}[(\hat{C}_n^{22} - C^{22})(\hat{C}_n^{22} -
C^{22})']\|_2^2 = O(1/nh^d)$;
\item[(c)] $\hat{C}_n \stackrel{p}{\rightarrow} C$.
\end{itemize}
\end{lemma}
\begin{pf}
This proof follows the techniques of \cite
{ruppertwand1994,bickelli2007}. We first prove part (c). Note that
\[
\hat{C}_n^{11} = \frac{1}{nh^{d-p}}\sum_{i=1}^n{K_h(x_i - x_0)},
\]
and consider its expectation
\begin{eqnarray*}
\mathbb{E}(\hat{C}_n^{11}) & = & \mathbb{E}\biggl(\frac
{1}{h^{d-p}}K_h(x_i -
x_0)\mathbh{1}\bigl(X \in(\mathcal{B}^p_{x,h^{1-\varepsilon}})\bigr
)\biggr) \\
&&{} + \mathbb{E}\biggl(\frac{h^p}{h^d}K_h(x_i - x_0)\mathbh{1}\bigl(X
\in(\mathcal
{B}^p_{x,h^{1-\varepsilon}})^c\bigr)\biggr) \\
& = & \int_{B^d_{0, h^{1-\varepsilon}}}{\frac{1}{h^d}K\biggl(\frac{\phi
(z) -
\phi(0)}{h}\biggr)F(z)\,dz} + o(h^{2+p}) \\
& = & \int_{\mathbb{R}^n}{K(d_u\phi\cdot u)F(0)\,du} + O(h^2),
\end{eqnarray*}
where we have used the assumption that $K(\cdot)$ is an even function,
$K'(\cdot)$ is an odd function, and $K''(\cdot)$ is an even function.

A similar calculation shows that for
\[
\hat{C}_n^{21} = \frac{1}{nh^{d+1-p}}\sum_{i=1}^n{K_h(x_i - x_0)(x_i - x_0)},
\]
we have that the expectation is
\[
\mathbb{E}(\hat{C}_n^{21}) = O(h) = o(1).
\]
And, a similar calculation shows that for
\[
\hat{C}_n^{22} = \frac{1}{nh^{d+2-p}}\sum_{i=1}^n{K_h(x_i - x_0)(x_i -
x_0)(x_i - x_0)'},
\]
we have that the expectation is
\[
\mathbb{E}(\hat{C}_n^{22}) = F(0)\,d_u\phi\cdot\biggl[\int_{\mathbb
{R}^d}{K\bigl(d_u\phi(0) \cdot u\bigr)uu'\,du}\biggr] \cdot d_u\phi' + O(h^2).
\]
The result in part (c) follows from the weak law of large numbers. The
last calculation also proves part (a).

Next, we prove part (b). For notational simplicity, let
\[
T_i = K_h(x_i - x_0)(x_i - x_0)(x_i - x_0)'.
\]
The variance is
\begin{eqnarray*}
\operatorname{Var}(\hat{C}_n^{22}) & = & \frac
{1}{n^2h^{4+2d-2p}}\operatorname{Tr}\bigl(n\bigl(\mathbb{E}(T_iT_i') -
\mathbb
{E}(T_i)\mathbb{E}(T_i)'\bigr) \\
&&\hspace*{70.1pt}{} + n(n-1)\bigl(\mathbb{E}(T_i T_j') - \mathbb
{E}(T_i)\mathbb{E}(T_j)'\bigr)\bigr).
\end{eqnarray*}
Since $T_i$ and $T_j$ are independent, it follows that $\mathbb
{E}(T_iT_j') - \mathbb{E}(T_i)\mathbb{E}(T_j)' = 0$. Next, note that
\begin{eqnarray*}
\mathbb{E}(T_i T_i') & = & h^{-2p}\biggl(\int_{B_{0,h^{1-\varepsilon
}}^d}\biggl(K\biggl(\frac{\phi(z) - \phi(0)}{h}\biggr)\biggr)^2[\phi(z)
- \phi(0)][\phi(z) - \phi(0)]' \\
&&\hspace*{64.2pt}{} \times[\phi(z) - \phi(0)][\phi(z) - \phi
(0)]'F(z)\,dz + o(h^{d+2})\biggr) \\
& = & h^{d+4-2p}\biggl(F(0)\int_{\mathbb{R}^d}\bigl(K(d_u\phi\cdot
u)\bigr)^2\,d_u\phi
\cdot uu' \cdot d_u\phi' \cdot
d_u\phi\\
&&\,\hspace*{123.5pt}{}\times uu'\cdot d_u\phi' \cdot du + O(h^2)\biggr).
\end{eqnarray*}
Thus, the variance is given by
\begin{eqnarray*}
\operatorname{Var}(\hat{C}_n^{22}) & = & \frac{1}{nh^{d}}\operatorname
{Tr}\biggl(F(0)\int_{\mathbb{R}^d}\bigl(K(d_u\phi\cdot u)\bigr)^2\,d_u\phi
\cdot uu'
\cdot d_u\phi' \cdot d_u\phi\\
&&\,\hspace*{119pt}{}\times uu'\cdot d_u\phi' \cdot du + o_p(1)\biggr).
\end{eqnarray*}
\upqed
\end{pf}
\begin{lemma}\label{theorem:proj} If the assumptions in Section \ref
{section:setup} hold, then the matrices $\hat{C}_n$, $C^{22}$, $\Pi$,
$\hat{\Pi}_n$, $\hat{P}_n$, and $P$ have the following properties:
\begin{enumerate}

\item[(a)] $\operatorname{rank}(C^{22}) = d$ and $\mathcal{R}(C^{22}) =
T_p\mathcal{M}$;

\item[(b)] $\mathcal{R}(\Pi) = \mathcal{N}(C^{22})$, $\mathcal{N}(\Pi) =
\mathcal{R}(C^{22})$ and $\mathcal{N}(\Pi) \cap\mathcal{N}(C^{22}) =
\{0\}$;

\item[(c)] $\|\hat{P}_n - P\|_2^2 = \|\hat{\Pi}_n - \Pi\|_2^2 = O_p(1/nh^d)$;

\item[(d)] $\mathbb{P}(\operatorname{rank}(\hat{C}_n + \lambda_n \hat
{P}_n/nh^{d+2}) = p+1) \rightarrow1$.
\end{enumerate}
\end{lemma}
\begin{pf}
To show property (a), we first show that for $M \in\mathbb{R}^{d
\times d}$, where
\[
M = \int_{\mathbb{R}^d}{K\bigl(d_u\phi(0) \cdot u\bigr)uu'\,du},
\]
we have that $\operatorname{rank}(M) = d$. To prove this, choose any $v
\in\mathbb{R}^d \setminus\{0\}$ and then consider the quantity
\[
v'Mv = \int_{\mathbb{R}^d}{K\bigl(d_u\phi(0) \cdot u\bigr)v'uu'v\,du}.
\]
By construction, $v'uu'v > 0$ almost everywhere. Additionally, since
$\phi$ is three times differentiable, we have that $K(d_u\phi(0) \cdot
u) > 0$ on a set of nonzero measure and $K(d_u\phi(0) \cdot u) \geq0$
elsewhere. Thus, $v'Mv > 0$ for all $v \in\mathbb{R}^d \setminus\{0\}
$. It follows that $M$ is symmetric and positive definite with
$\operatorname{rank}(M) = d$. Since $\mathcal{M}$ is a $d$-dimensional
manifold, we have that $\operatorname{rank}(d_u\phi) = d$ by Corollary
8.4 of \cite{lee2003}. The Sylvester Inequality \cite{sastry1999}
implies that
\[
\operatorname{rank}(C_{22}) = \operatorname{rank}(d_u\phi M d_u\phi') = d,
\]
and this implies that
\[
\mathcal{R}(C_{22}) = \mathcal{R}(d_u\phi M d_u\phi') = \mathcal{R}(d_u
\phi).
\]
However, $\mathcal{R}(d_u \phi) = T_p\mathcal{M}$, where we take $p =
x_0$. This proves the result.

We next consider property (b). We have that
\[
\sigma_1, \ldots, \sigma_d \neq0 \quad\mbox{and}\quad\sigma_{d+1} =
\cdots=
\sigma_p = 0,
\]
because $\operatorname{rank}(C_{22}) = d$ by property (a). Thus, the
null-space of $C_{22}$ is given by the column-span of $U_N$; however,
the construction of $P$ implies that the column-span of $U_N$ is the
range-space of $P$. Ergo, $\mathcal{R}(P) = \mathcal{N}(C_{22})$. Note
that the column-span of $U_R$ belongs to the null-space of $P$, because
each column in $U_R$ is orthogonal---by property of the SVD---to each
column in $U_N$. Thus, we have the dual result that $\mathcal{N}(P) =
\mathcal{R}(C_{22})$. The orthogonality of $U_R$ and $U_N$ due to the
SVD implies that $\mathcal{N}(P) \cap\mathcal{N}(C_{22}) = \{0\}$.

Now, we turn to property (c). For $h = \kappa n^{-1/(d+4)}$, Lemma \ref
{theorem:gram} says that $\|\hat{C}_n^{22} - C_{22}\|_F^2 =
O_p(1/nh^d)$. The result follows from Corollary 3 of \cite
{elkaroui2007}, by the fact that $\mathbb{I} - P_X$ is the projection
matrix onto the null-space of $X$, and by the equivalence:
\[
\|P_X - P_Z\|_2^2 \equiv\|{\sin\Theta}[\mathcal{R}(X), \mathcal
{R}(Z)]\|,
\]
where $P_X$ is a projection matrix onto the range space of $X$ \cite
{stewartsun1990}.

Lastly, we deal with property (d). Lemma \ref{theorem:gram} shows that
\[
C_n \stackrel{p}{\rightarrow} C = \left[\matrix{ C_{11} & 0 \cr0 &
C_{22}}\right].
\]
Since $F(0) \neq0$ by assumption, $C_{11} \neq0$; thus, $\operatorname
{rank}(C) = 1 + \operatorname{rank}(C_{22})$. Since $\mathcal{N}(P)
\cap\mathcal{N}(C_{22}) = \{0\}$, we have that $\operatorname
{rank}(C_{22} + \lambda_nP/nh^{d+2}) = p$. Consequently, $\operatorname
{rank}(C) = p+1$. Next, consider the expression
\begin{eqnarray*}
&& \|C_n + \lambda_n\tilde{P}_n/nh^{d+2} - C - \lambda_n\tilde
{P}/nh^{d+2}\|_2^2 \\
&&\qquad\leq\|\hat{C}_n - C\|_2^2 + \frac{\lambda_n}{nh^{d+2}}\|\tilde{P}_n
- \tilde{P}\|_2^2\\
&&\qquad\leq O_p(h^2) + O_p\biggl(\frac{\lambda_n}{n^2h^{2d+2}}\biggr)
\\
&&\qquad\leq o_p(h).
\end{eqnarray*}
Weyl's theorem \cite{bhatia2007} implies that
%
%
\begin{equation}
\label{eqn:weyl}
\|\sigma_i(C_n + \lambda_n\tilde{P}_n/nh^{d+2}) - \sigma_i(C + \lambda
_n\tilde{P}/nh^{d+2})\|_2^2 \leq o_p(h).
\end{equation}
Note that $\sigma_i(C + \lambda_n\tilde{P}/nh^{d+2})$ is nondecreasing
because $\lambda_n/nh^{d+2}$ is nondecreasing. Define
\[
\eta= \min\bigl(\sigma_i\bigl(C + \lambda_n\tilde{P}/n^{2/(d+4)}\bigr
)\bigr),
\]
and consider the probability
%
%
\begin{eqnarray}\label{eqn:bound}
&& \mathbb{P}\bigl(\operatorname{rank}(C_n + \lambda_n \tilde{P}_n/nh^{d+2})
= p+1\bigr) \nonumber\\
&&\qquad\geq\mathbb{P}\bigl(\bigl|\sigma_i(C_n + \lambda_n\tilde{P}_n/nh^{d+2})
\nonumber\\
&&\qquad\quad\hspace*{11.6pt}{} -\sigma_i(C + \lambda_n\tilde
{P}/nh^{d+2})\bigr| \leq\eta, \forall i\bigr)
\\
&&\qquad\geq\sum_{i=1}^{p+1}\mathbb{P}\bigl(\bigl|\sigma_i(C_n +
\lambda_n\tilde
{P}_n/nh^{d+2}) \nonumber\\
&&\qquad\quad\hspace*{29.2pt}{} -\sigma_i(C + \lambda_n\tilde
{P}/nh^{d+2})\bigr| \leq\eta\bigr) - p. \nonumber
\end{eqnarray}
The result follows from (\ref{eqn:weyl}) and (\ref{eqn:bound}).
\end{pf}

For notational convenience, we define
\[
B_n = h^p\bigl(f(X) - \beta X_x\bigr)'W_{x_0}X_xH^{-1/2}
\]
and
\[
M = \tfrac{1}{2}[\partial_i\partial_j(f \circ\phi)- d_xf \cdot\partial
_i\partial_j\phi].
\]
Then, we have the following result concerning the asymptotic bias of
the estimator:
\begin{lemma}
If $h = \kappa n^{-1/(d+4)}$, then $B_n \stackrel{p}{\rightarrow} B$, where
%
%
\begin{equation}
\label{eqn:bee}
B = \biggl[\matrix{ \displaystyle\kappa F(0)\int_{\mathbb
{R}^d}{K(d_u\phi\cdot u)uu'\,du}\,\displaystyle M'
& \mathbb{O}_{p \times1} }\biggr].
\end{equation}
%
\end{lemma}
\begin{pf}
First, recall the Taylor polynomial of the pullback of $f$ to $z$:
\[
f(\phi(z)) = f(\phi(0)) + d_xf \cdot d_u\phi\cdot z + \tfrac
{1}{2}\,\partial_i\partial_j(f \circ\phi) \cdot zz' + o(\|z\|^2),
\]
where we have performed a pullback of $d_xf$ from $T_x^*\mathcal{M}$ to
$T_x^*\mathcal{U}$. In the following expression, we set $z = hu$:
\begin{eqnarray*}
&& f(\phi(z)) - f(x_0) - d_xf \cdot[\phi(z) - \phi(0)] \\
&&\qquad= \frac{h^2}{2}[\partial_i\partial_j(f \circ\phi)- d_xf \cdot
\partial_i\partial_j\phi]uu' + o(\|hu\|^2) \\
&&\qquad= h^2Buu' + o(\|hu\|^2).
\end{eqnarray*}
Because $\beta= [ f(x_0) \enskip d_xf ]'$, we can rewrite the
expectation of the expression as
\begin{eqnarray*}
\mathbb{E}(B_n) & = & \mathbb{E}\left(h^{p-d} K_H(x - x_0)\bigl(f(x) -
x_{x_0}\beta
\bigr)'x_{x_0}\left[\matrix{ 1 & 0 \cr0 & 1/h^2 \mathbb{I}}\right]\right
)H^{1/2} \\
& = & \biggl(\int_{\mathbb{R}^d}\biggl\{K(d_u\phi\cdot u)h^2u'uM' \biggl
[\matrix{ 1 &
\displaystyle\frac{1}{h}d_u\phi\cdot u + \frac{1}{2}\,\partial_i\partial
_j\phi\cdot
uu' }\biggr] \\
&&\,\hspace*{56.9pt}\hspace*{89.2pt}{} \times\bigl(F(0) + hd_uF(0) \cdot u\bigr)\biggr\}\,
du\\
&&\,\hspace*{149pt}\hspace*{89.2pt}{} + o(h^2)\biggr)H^{1/2} \\
& = &\sqrt{nh^d}h^2\biggl[\matrix{ \displaystyle F(0)\int_{\mathbb
{R}^d}{K(d_u\phi\cdot
u)uu'\,du}M' + o(1) & O\bigl(\sqrt{h^2}\bigr) }\biggr],
\end{eqnarray*}
where the last line follows because of the odd symmetries in the
integrand. Since $h = \kappa n^{-1/(d+4)}$, this expectation becomes
\[
\mathbb{E}(B_n) = B + o(1)\mathbh{1}_{1 \times(p+1)}.
\]
The result follows from application of the weak law of large numbers.
\end{pf}

Let
%
%
\begin{equation}
\label{eqn:vee}
V = F(0)\int_{\mathbb{R}^d}{\bigl(K(d_u\phi\cdot u)\bigr)^2\left[\matrix
{ 1 & 0 \cr0
& d_u\phi\cdot uu' \cdot d_u\phi}\right]\,du},
\end{equation}
then the following lemma describes the asymptotic distribution of the
error residuals.
\begin{lemma}
If $h = \kappa n^{-1/(d+4)}$, then
\[
h^p\varepsilon'W_{x_0}X_xH^{-1/2} \stackrel{d}{\rightarrow} \mathcal
{N}(0, \sigma^2V).
\]
\end{lemma}
\begin{pf}
Since $\mathbb{E}(\varepsilon) = 0$ and $\varepsilon$ is independent of
$x$, we have that
\[
\mathbb{E}\bigl(\sqrt{n}h^p\varepsilon K_H(x - x_0)x_{x_0}H^{-1/2}\bigr
) = 0.
\]
The variance of this quantity is
\begin{eqnarray*}
&& \operatorname{Var}\bigl(h^p\sqrt{n}\varepsilon K_H(x -
x_0)x_{x_0}H^{-1/2}\bigr) \\
&&\qquad= \mathbb{E}\bigl(\bigl(h^p\sqrt{n}\varepsilon K_H(x -
x_0)x_{x_0}H^{-1/2}\bigr)'\bigl(h^p\sqrt{n}\varepsilon K_H(x -
x_0)x_{x_0}H^{-1/2}\bigr)\bigr) \\
&&\qquad= nh^{2p}\sigma^2\mathbb{E}\bigl(\bigl(K_H(x - x_0)\bigr
)^2\left[\matrix{ 1 & (x - x_0)'
}\right]'H^{-1}\left[\matrix{ 1 & (x - x_0)' }\right]\bigr) \\
&&\qquad= \sigma^2\left\{\int_{\mathbb{R}^d}\bigl(K(d_u\phi\cdot
u)\bigr)^2\left[\matrix{ 1 &
\biggl(d_u\phi\cdot u + \dfrac{h}{2}\,\partial_i\partial_j \phi\cdot
uu'\biggr)'
\cr
\cdot& d_u\phi\cdot uu' \cdot d_u\phi}\right]\right. \\
&&\,\qquad\quad\hspace*{87.7pt}{} \times\left.\bigl(F(0)+ hd_uF(0) \cdot
u\bigr)\,du + o(h^2)
\vphantom{\matrix{ 1 &
\biggl(d_u\phi\cdot u + \dfrac{h}{2}\,\partial_i\partial_j \phi\cdot
uu'\biggr)'
\cr
\cdot& d_u\phi\cdot uu' \cdot d_u\phi}}\right\} \\
&&\qquad= \sigma^2\bigl(V + o(h)\mathbb{I}\bigr).
\end{eqnarray*}
Thus, the central limit theorem implies that
\[
\frac{h^p\sqrt{n}\varepsilon'W_{x_0}X_xH^{-1/2}}{\sqrt{n}} \stackrel
{d}{\rightarrow} \mathcal{N}(0,\sigma^2V).
\]
%
\upqed
\end{pf}
\begin{pf*}{Proof of Theorem \ref{theorem:nede}}
This proof follows the framework of \cite{knightfu2000,zou2006} but
with significant modifications to deal with our estimator. For
notational convenience, we define the indices
of $\beta$ such that: $\beta_0 = f(x_0)$ and $[\beta_1 \enskip\cdots
\enskip\beta_p ] = d_xf$.
Let $\tilde{\beta} = \beta+ H^{-1/2}u$ and
\begin{eqnarray*}
\Psi_n(u) &=& h^p\bigl\|W_{x_0}^{1/2}\bigl(Y - X_{x_0}(\beta+
H^{-1/2}u)\bigr)\bigr\|_2^2 \\
&&{} + \lambda_n\|P_n \cdot(\beta+ H^{-1/2}u)\|_2^2.
\end{eqnarray*}
Let $\hat{u}^{(n)} = \arg\min\Psi_n(u)$; then $\hat{\beta}^{(n)} =
\beta+ H^{-1/2}\hat{u}^{(n)}$. Note that $\Psi_n(u) - \Psi_n(0) =
V_4^{(n)}(u)$, where
\begin{eqnarray*}
V_4^{(n)}(u) &=& u'H^{-1/2}(h^pX_{x_0}'W_{x_0}X_{x_0} + \lambda_n
P_n)H^{-1/2}u \\
&&{}+ 2\bigl(h^p(Y - X_{x_0}\beta)'W_{x_0}X_{x_0} + \lambda_n\beta
'P_n\bigr)H^{-1/2}u.
\end{eqnarray*}
If $\lambda_n/nh^{d+2} \rightarrow\infty$ and $h\lambda_n/nh^{d+2}
\rightarrow0$, then for every $u$
\[
\lambda_n\beta'P_nH^{-1/2}u = \lambda_n\beta'Pu/\sqrt{nh^{d+2}}O_p(1) +
\lambda_nh/nh^{d+2}O_p(1),
\]
where we have used Lemma \ref{theorem:proj}. It follows from the
definition of $\beta$ (\ref{eqn:beta}) and Lemma \ref{theorem:proj}
that $\beta'P \equiv0$; thus, $\lambda_n\beta'P_nH^{-1/2}u = \lambda
_nh/nh^{d+2}O_p(1) = o_p(1)$.
For all $u \in\mathcal{N}(P)$, we have
\[
\lambda_n/nh^{d+2}u'P_nu = \lambda_n/nh^{d+2}O_p(1/nh^d) = o_P(h\lambda
_n/nh^{d+2}),
\]
and for all $u \notin\mathcal{N}(P)$, we have
\[
\lambda_n/nh^{d+2}u'P_nu = \lambda_n/nh^{d+2}u'P_nuO_p(1) \rightarrow
\infty.
\]

Let $W \sim\mathcal{N}(0, \sigma^2V)$. Then, by Slutsky's theorem we
must have that $V_4^{(n)}(u) \stackrel{d}{\rightarrow} V_4(u)$ for
every $u$, where
\[
V_4(u) = \cases{
u'Cu - 2u'(W + B), &\quad if $u \in\mathcal{N}(P)$, \cr
\infty, &\quad otherwise.}
\]
Lemma 5 shows that $V_4^{(n)}(u)$ is convex with high-probability, and
Lemma 5 also shows that $V_4(u)$ is convex. Consequently, the unique
minimum of $V_4(u)$ is given by $u = C^{\dag}(W+B)$, where $C^{\dag}$
denotes the Moore--Penrose pseudoinverse of $C$. Following the
epi-convergence results of \cite{geyer1994,knightfu2000}, we have that
$\hat{u}^{(n)} \stackrel{d}{\rightarrow} C^{\dag}(W+B)$. This proves
asymptotic normality of the estimator, as well as convergence in probability.

The proof for the NALEDE estimator comes for free. The proof
formulation that we have used for the consistency of nonparametric
regression in (\ref{eqn:nede}) allows us to trivially extend the proof
of \cite{zou2006} to prove asymptotic normality and consistency.
\end{pf*}
\begin{lemma}
\label{lemma:invmat}
Consider $A_n, B_n \in\mathbb{R}^{p_n \times p_n}$ that are symmetric,
invertible matrices. If $\|A_n - B_n\|_2 = O_p(\gamma_n)$, $\|A_n^{-1}\|
_2 = O_p(1)$ and $\|B_n^{-1}\|_2 = O_p(1)$, then $\|A_n^{-1} -
B_n^{-1}\|_2 = O_p(\gamma_n)$.
\end{lemma}
\begin{pf}
Consider the expression
\begin{eqnarray*}
\|A_n^{-1} - B_n^{-1}\|_2 & = & \|A_n^{-1}(B_n - A_n)B_n^{-1}\|_2 \\
& \leq& \|A_n^{-1}\|_2 \cdot\|A_n - B_n\|_2 \cdot\|B_n^{-1}\|_2,
\end{eqnarray*}
where the last line follows because the induced, matrix norm \mbox{$\|
\cdot\|
_2$} is sub-multiplicative for square matrices.
\end{pf}
%
%
%
\begin{pf*}{Proof of Theorem \ref{theorem:ede}}
Under our set of assumptions, the results from \cite{bickellevina2008} apply:
%
%
\begin{eqnarray}
\|T_t(X'X/n) - (\Sigma_{\xi} + \sigma_{\nu}^2\mathbb{I})\|_2 & = &
O_p\Biggl(c_n\sqrt{\frac{\log p}{n}}\Biggr), \\
\label{eqn:2ndone}
\|T_t(X'Y/n) - \Sigma_{\xi}\overline{\beta}\|_2 & = &
O_p\Biggl(c_n\sqrt{\frac{\log p}{n}}\Biggr).
\end{eqnarray}
An argument similar to that given in Lemma \ref{theorem:proj} implies that
\[
\|\hat{P}_n - P_n\| = O_p\bigl(c_n\sqrt{\log p/n}\bigr).
\]
Consequently, it holds that
%
%
\begin{eqnarray}
\label{eqn:3rdone}
&&
\|T_t(X'X/n) - \sigma_{\nu}^2\mathbb{I} + \lambda_n \hat{P}_n - (\Sigma
_{\xi} + \lambda_n P_n)\|_2 \nonumber\\[-8pt]\\[-8pt]
&&\qquad= O_p\Biggl(c_n(\lambda_n+1)\sqrt{\frac{\log p}{n}}\Biggr
).\nonumber
\end{eqnarray}

Next, observe that
\[
\Sigma_{\xi} + \lambda_n P_n = \left[\matrix{ U_R & U_N }\right]
\operatorname
{diag}(\sigma_1, \ldots, \sigma_d, \lambda_n, \ldots, \lambda_n)\left
[\matrix
{ U_R & U_N }\right]'.
\]
Recall that we only consider the case in which $d < p$. We have that:
\begin{itemize}
\item[(a)] $\|(\Sigma_{\xi} + \lambda_n P_n)^{-1}\|_2 = O(1)$, because
of (\ref{eqn:sigbound});
\item[(b)] $\|\Sigma_{\xi}^{\dag} - (\Sigma_{\xi} + \lambda_n
P_n)^{-1}\|
_2 = O_p(1/\lambda_n)$.
\end{itemize}
Weyl's theorem \cite{bhatia2007} and (\ref{eqn:3rdone}) imply that
\[
\bigl\|\bigl(T_t(X'X/n) - \sigma_{\nu}^2\mathbb{I} + \lambda_n \hat
{P}_n\bigr)^{-1}\bigr\|
_2 = O_p(1).
\]
Additionally, Lemma \ref{lemma:invmat} implies that
\[
\bigl\|\bigl(T_t(X'X/n) - \sigma_{\nu}^2\mathbb{I} + \lambda_n \hat
{P}_n\bigr)^{-1} -
(\Sigma_{\xi} + \lambda_n P_n)^{-1}\bigr\| = O_p\Biggl(c_n\lambda_n\sqrt
{\frac
{\log p}{n}}\Biggr).
\]

Note that the solution to the estimator defined in (\ref{eqn:ede}) is:
\[
\hat{\beta} = \bigl(T_t(X'X/n) - \sigma_{\nu}^2\mathbb{I} + \lambda_n
\hat
{P}_n\bigr)^{-1}T_t(X'Y/n).
\]
Next, we define
\[
\beta^{(n)} \triangleq\bigl(T_t(X'X/n) - \sigma_{\nu}^2\mathbb{I} +
\lambda
_n \hat{P}_n\bigr)^{-1}\Sigma_{\xi}\overline{\beta},
\]
and note that the projection matrix onto the range space of $\Sigma_{\xi
}$ is given by $P_{\Sigma_{\xi}} = \Sigma_{\xi}^{\dag}\Sigma_{\xi}$.
Thus, $\beta= P_{\Sigma_{\xi}}\overline{\beta} = \Sigma_{\xi}^{\dag
}\Sigma_{\xi}\overline{\beta}$. Consequently, we have that
%
%
\begin{eqnarray} \label{eqn:betasuff}\qquad
\|\hat{\beta} - \beta\|_2
& \leq&\bigl\|\hat{\beta} - \beta^{(n)}\bigr\|_2 + \bigl\|\beta^{(n)}
- \beta\bigr\|_2
\nonumber\\
& \leq&\bigl\|\bigl(T_t(X'X/n) - \sigma_{\nu}^2\mathbb{I} + \lambda_n
\hat
{P}_n\bigr)^{-1}\bigr\|_2\cdot\|T_t(X'Y/n) -
\Sigma_{\xi}\overline{\beta}\|_2\\
&&{} + \bigl\|\bigl(T_t(X'X/n) - \sigma_{\nu
}^2\mathbb{I} + \lambda_n \hat{P}_n\bigr)^{-1}
- \Sigma_{\xi}^{\dag}\bigr\|_2 \cdot\|\Sigma_{\xi}\overline{\beta}\|
_2,\nonumber
\end{eqnarray}
where the inequality comes about because $\|\cdot\|_2$ is an induced,
matrix norm and the expressions are of the form $\mathbb{R}^{p \times
p}(\mathbb{R}^{p \times p}\mathbb{R}^p)$. Recall that for symmetric
matrices, $\|A\|_1 = \|A\|_\infty$; ergo, $\|A\|_2 \leq\sqrt{\|A\|_1\|
A\|_\infty} = \|A\|_1$. Because of (\ref{eqn:sparselinear}), we can use
this relationship on the norms to calculate that $\|\Sigma_{\xi}\|_2 =
O(c_n)$ and $\|\overline{\beta}\| = O(c_n)$. Consequently,
\[
\mbox{(\ref{eqn:betasuff})} \leq O_p\Biggl(c_n\lambda_n\sqrt{\frac{\log
p}{n}}\Biggr) +
O_p(c_n^2/\lambda_n).
\]
The result follows from the relationship
\[
\lambda_n = O\biggl(\sqrt{c_n}\biggl(\frac{n}{\log p}\biggr)^{1/4}\biggr).
\]
\upqed
\end{pf*}

We can show that the bias of the terms of the nonparametric exterior
derivative estimation goes to zero at a certain rate.
\begin{lemma}
\label{theorem:crate}
Under the assumptions of Section \ref{section:lpsnnpc}, we have that
\begin{eqnarray*}
|\mathbb{E}([\hat{C}_n]_{ij}) - [C_n]_{ij}| & = & O(h^2c_n^2(2\Omega
)^{2d}), \\
|\mathbb{E}([\hat{R}_n]_{ij}) - [R_n]_{ij}| & = & O(h^2c_n^2(2\Omega)^{2d}),
\end{eqnarray*}
where $i, j$ denote the components of the matrices. Similarly, we have that
\begin{eqnarray*}
\operatorname{Var}([n\hat{C}_n]_{ij}) & = & O(1/h^d), \\
\operatorname{Var}([n\hat{R}_n]_{ij}) & = & O(1/h^d).
\end{eqnarray*}
%
\end{lemma}
\begin{pf}
By the triangle inequality and a change of variables,
\begin{eqnarray*}
\operatorname{Bias}(\hat{C}_n^{11}) & = &
\biggl|\int_{\mathcal{B}^d_{0, \Omega/h}} \frac{1}{h^d}K\biggl(\frac
{\phi(z) -
\phi(0)}{h}\biggr)F(z)\,dz - \int_{\mathcal{B}^d_{0, \Omega}} K(d_u\phi
\cdot
u)F(0)\,du\biggr| \\
& \leq& \biggl|\int_{\mathcal{B}^d_{0, \Omega}} \biggl[K\biggl(\frac
{\phi(hu) - \phi
(0)}{h}\biggr) - K(d_u \phi\cdot u)\biggr]F(hu)\,du\biggr| \\
&&{} + \biggl|\int_{\mathcal{B}^d_{0, \Omega}} K(d_u \phi\cdot u)[F(hu) -
F(0)]\,du\biggr| = T_1 + T_2.
\end{eqnarray*}
The Taylor remainder theorem implies that
\begin{eqnarray*}
&&K\biggl(\frac{\phi(hu) - \phi(0)}{h}\biggr)
\\
&&\qquad= K(d_u \phi\cdot u)\\
&&\qquad\quad{} + \partial_k
K(d_u\phi\cdot u) \times(h\,\partial_{ij}\phi^k|_0 u^i u^j/2 + h^2\,\partial
_{ijm}\phi
^k|_w u^i u^j u^m/6) \\
&&\qquad\quad{} + \partial_{kl}K(v)/2 \times(h\,\partial_{ij}\phi^k|_0
u^i u^j/2
+ h^2\,\partial_{ijm}\phi^k|_w u^i u^j u^m/6)\\
&&\qquad\quad\hspace*{9.9pt}{}\times(h\,\partial_{ij}\phi
^l|_0 u^i u^j/2 + h^2\,\partial_{ijm}\phi^l|_w u^i u^j u^m/6),
\end{eqnarray*}
where $w \in\overline{\mathcal{B}^d_{0, \Omega}}$ and $v \in\overline
{\mathcal{B}^d_{d_u\phi\cdot u, h\,\partial_{ij}\phi^k|_0 u^i u^j/2 +
h^2\,\partial_{ijm}\phi^k|_w u^i u^j u^m/6}}$, and
\[
F(hu) = F(0) + h\,\partial_iF|_0u^i + h^2\,\partial_{ij}F|_vu^iu^j/2,
\]
where $v \in(0, hu)$.

The odd-symmetry components of the integrands of $T_1$ and $T_2$ will
be equal to zero, and so we only need to consider even-symmetry terms
of the integrands. Recall that $K(\cdot), \partial_k K(\cdot), \partial
_{kl} K(\cdot)$ are, respectively, even, odd and even. By the sparsity
assumptions, we have that
\begin{eqnarray*}
T_1 & = & O(h^2d^6c_n^2(2\Omega)^d), \\
T_2 & = & O(h^2d^2(2\Omega)^d).
\end{eqnarray*}
Consequently, $T_1 + T_2 = O(h^2d^6c_n^2(2\Omega)^{d}) =
O(h^2c_n^2(2\Omega)^{2d})$.

We can compute the variance of $n\hat{C}_n^{11}$ to be
\begin{eqnarray*}
\operatorname{Var}(n\hat{C}_n^{11}) & = &\int_{\mathcal{B}^d_{0, \Omega
/h}}h^{-2d}\bigl[K\bigl(\bigl(\phi(z) - \phi(0)\bigr)/h\bigr)\bigr
]^2(F(z))^2\, dz \\
&&{} - (\mathbb{E}(n\hat{C}_n^{11}))^2 \\
& = &h^{-d}\int_{\mathcal{B}^d_{0, \Omega}}\bigl[K\bigl(\bigl(\phi(hu)
- \phi
(0)\bigr)/h\bigr)\bigr]^2(F(hu))^2 \,dy \\
&&{} - (\mathbb{E}(n\hat{C}_n^{11}))^2 \\
& = &O(1/h^d).
\end{eqnarray*}

The remainder of the results follow by similar, lengthy calculations.
Note that for the variance of terms involving $Y_i$, a $\sigma^2$
coefficient appears, but this is just a finite-scaling factor which is
irrelevant in $O$-notation.
\end{pf}
\begin{pf*}{Proof of Theorem \ref{theorem:lpsnthresh}}
The key to this proof is to provide an exponential concentration
inequality for the terms in $\hat{C}_n$ and $\hat{R}_n$. Having done
this, we can then piggyback off of the proof in \cite{bickellevina2008}
to immediately get the result. The proofs\vspace*{2pt} for $\hat{C}_n$
and $\hat
{R}_n$ are identical; so we only do the proof for $\hat{C}_n$.

Using the Bernstein inequality \cite{lugosi2006} and the union bound,
\begin{eqnarray*}
&&\mathbb{P}\Bigl({\max_{i,j}}\|[\hat{C}_n]_{ij} - \mathbb{E}[\hat
{C_n}]_{ij}\|
> t \Bigr)\\
&&\qquad\leq2p^2\exp\biggl(-\frac{nt^2}{2\operatorname{Var}(n[\hat
{C}_n]_{ij}) +
\max(|n[\hat{C}_n]_{ij}|) 2t/3}\biggr).
\end{eqnarray*}
Since the $i$th component of $X$ obeys: $|[X]_i| \leq M$, it follows that
\[
\max(|n[\hat{C}_n]_{ij}|) = 2M/h^\eta,
\]
where $\eta\in\{0, 1, 2\}$ depending on $i$ and $j$. Using this bound
and Lemma \ref{theorem:crate} gives
\[
{\max_{i,j}}|[\hat{C}_n]_{ij} - \mathbb{E}[\hat{C_n}]_{ij}| = O_p\bigl
(\sqrt
{\log p/nh^d}\bigr).
\]
Recall that
\[
{\max_{i,j}}|[\hat{C}_n]_{ij} - [C_n]_{ij}| \leq{\max_{i,j}}|[\hat
{C}_n]_{ij} - \mathbb{E}[\hat{C_n}]_{ij}| + {\max_{i,j}}|\mathbb
{E}([\hat
{C}_n]_{ij}) - [C_n]_{ij}|.
\]
However, this second term is $o(\sqrt{\log p/nh^d})$. Consequently,
%
%
\begin{equation}
\label{eqn:rate}
{\max_{i,j}}|[\hat{C}_n]_{ij} - [C_n]_{ij}| = O_p\bigl(\sqrt{\log
p/nh^d}\bigr).
\end{equation}
Using (\ref{eqn:rate}), we can follow the proof of Theorem 1 in \cite
{bickellevina2008} to prove the result.
\end{pf*}
\end{appendix}


%
\printaddresses

\end{document}